\RequirePackage{fix-cm}
\documentclass[9pt]{article}       
\usepackage{algorithm}
\usepackage{algorithmic}
\usepackage{hyperref}

\usepackage{amsmath}
\usepackage{amssymb}
\usepackage{tikz}
\usepackage{epstopdf}
\usetikzlibrary{arrows,matrix,positioning,backgrounds,trees}
\usepackage{subcaption}
\usepackage{graphicx}

\newcommand{\norm}[1]{||#1||}
\newcommand{\bnorm}[1]{\left|\left|#1\right|\right|}
\newcommand{\trace}[1]{\text{Tr}\left(#1\right)}
\newcommand{\vect}[1]{\text{vec}(#1)}
\newcommand{\innerp}[1]{\langle #1\rangle}

\DeclareMathOperator*{\argmin}{\arg\!\min}

\newcommand{\mbf}[1]{\mathbf{#1}}
\newcommand{\mcal}[1]{\mathcal{#1}}
\newcommand{\R}{\mathbb R}
\newcommand{\st}{\text{s.t.}}
\newcommand{\Perm}{\text{Perm}(n)}

\newcommand{\etal}{\emph{et al}}

\newtheorem{problem}{Problem}

\providecommand{\keywords}[1]{\textbf{Keywords:} #1}

\begin{document}

\title{Semidefinite Programming Approach for the Quadratic Assignment Problem with a Sparse Graph}

\author{Jos\'e F. S. Bravo Ferreira \thanks{Program in Applied and Computational Mathematics, Princeton University, Princeton, NJ 08544, USA (\texttt{josesf@princeton.edu})}	
	 \and Yuehaw Khoo \thanks{Department of Mathematics, Stanford University, Stanford, CA 94305, USA (\texttt{ykhoo@stanford.edu})}	
	 \and Amit Singer \thanks{Department of Mathematics and Program in Applied and Computational Mathematics, Princeton University, Princeton, NJ 08544, USA (\texttt{josesf@math.princeton.edu})}}

%


\maketitle

\begin{abstract}
The matching problem between two adjacency matrices can be formulated as the NP-hard quadratic assignment problem (QAP). Previous work on semidefinite programming (SDP) relaxations to the QAP have produced solutions that are often tight in practice, but such SDPs typically scale badly, involving matrix variables of dimension $n^2$ where $n$ is the number of nodes. To achieve a speed up, we propose a further relaxation of the SDP involving a number of positive semidefinite matrices of dimension $\mcal{O}(n)$ no greater than the number of edges in one of the graphs. The relaxation can be further strengthened by considering cliques in the graph, instead of edges. The dual problem of this novel relaxation has a natural three-block structure that can be solved via a convergent Augmented Direction Method of Multipliers (ADMM) in a distributed manner, where the most expensive step per iteration is computing the eigendecomposition of matrices of dimension $\mcal{O}(n)$. The new SDP relaxation produces strong bounds on quadratic assignment problems where one of the graphs is sparse with reduced computational complexity and running times, and can be used in the context of nuclear magnetic resonance spectroscopy (NMR) to tackle the assignment problem.\\

\end{abstract}

\section{Introduction}\label{sec:introduction}
Given two graphs, $\mathcal{G}_A$ and $\mathcal{G}_B$, with adjacency matrices $A$ and $B\in \R^{n\times n}$, respectively, the graph matching problem is that of finding a permutation matrix $P\in\Perm$ such that $A=PBP^T$ or $AP = PB$. This problem is also known as the graph isomorphism problem. While a breakthrough result in \cite{GIquasipolynomial} shows that the graph isomorphism problem has a worst-case time complexity of $\exp((\log n)^{\mcal{O}(1)})$, for most practical situations this problem can be solved very efficiently.

The problem becomes more complicated if $A$ and $B$ cannot be matched exactly. In this case, one needs to find the permutation matrix such that
$$\norm{AP-PB}$$
is minimized, where the typical choices of the norm $\norm{\cdot}$ are the entry-wise $\ell_1$ or $\ell_\infty$ norms and the Frobenius norm. Since the domain of the problem $\Perm$ is non-convex and combinatorially large, convex relaxation methods have been applied to search for the global optimum efficiently. In \cite{ramana1994fractional} and \cite{almohamad1993linear}, a convex relaxation is derived by relaxing the set of permutation matrices $\Perm$ to its convex hull, i.e. the set of doubly stochastic matrices. Under the $\ell_1$ or $\ell_\infty$ norms, the relaxed problem is a linear program (LP), while under the Frobenius norm it is a quadratic program (QP) with linear constraints.
In \cite{Aflalo}, the authors proved that this relaxation exactly solves the original problem if the graphs are isomorphic and friendly \footnote{A graph is called \emph{friendly} if its adjacency matrix has a simple spectrum and eigenvectors orthogonal to $\mbf 1_n$ \cite{Aflalo}.}. Furthermore, it produces an approximate isomorphism in the case of inexact matching of strongly friendly graphs. However, the need for a friendly graph is a rather strong condition. In particular, it is proven that this type of relaxation almost always fails to find the correct permutation for certain correlated Bernoulli random graphs, even for the case of exact graph matching \cite{Lyzinski}. In practice, we have found that this relaxation quickly loses its tightness in the presence of outlier-type noise.

On the other hand, the graph matching problem
\begin{align*}
\arg\min_P \norm{AP-PB}_F^2 &= \arg\min_P \trace{P^TA^TAP - 2PBP^TA^T + PBB^TP^T} \\
&=\arg\max_P \trace{PBP^TA^T}
\end{align*}
can be viewed as a special case of the Quadratic Assignment Problem (QAP). The QAP was first presented in \cite{Koopmans} and is known to be NP-hard (further, the $\epsilon$-approximation problem is also NP-hard \cite{Sahni}). It encodes a number of interesting problems, such as the traveling salesman problem (TSP) and the max clique problem (see, for example, \cite{Loiola} for a review of the applications of the QAP). The quadratic nature of the QAP invites a number of proposals to use semidefinite programming relaxations to attack the problem. The seminal SDP relaxation in \cite{Wolkowicz} has proven remarkably tight by achieving the optimal solution in several problem instances in the QAP library (QAPLIB, \cite{Burkard}). In \cite{Kovalsky}, the authors describe a very similar relaxation to tackle a shape matching problem in computer graphics. However, this convex relaxation introduces a semidefinite matrix variable of size $n^2\times n^2$, greatly hindering its use in practice. More recently, the alternating direction method of multipliers (ADMM) has been applied to ease the computational burden of solving this SDP, allowing problems with $n=30$ to be solved in a few minutes \cite{QAPadmmSDP}, but it remains a challenging problem to tackle, as it requires an eigendecomposition of a very large matrix at each iteration.

To make solving the QAP using SDP feasible, a few other convex relaxations have been proposed where the PSD variables have size $\mcal{O}(n)\times \mcal{O}(n)$. In particular, the relaxations in \cite{Peng1} and \cite{Peng2} achieve this by splitting $B$ into the difference of two PSD matrices, while those in \cite{deKlerk1} and \cite{deKlerk2} use the spectral decomposition of $PBP^T$. In the latter work, the symmetry of the matrix $A$ under graph automorphism can be used to achieve a significant reduction in problem size, and in special problem instances, such as the TSP, which possesses a cyclic symmetry, the SDP can be reduced to a linear program.

In section \ref{sec:related} we summarize the related works, highlighting the relaxation in \cite{Wolkowicz} that we use as a foundation for a novel edge-based convex relaxation. Section \ref{sec:csdp} presents the details of this new relaxation, extending it beyond edges to arbitrarily-sized cliques, and section \ref{sec:admm} demonstrates how one can use ADMM to solve the dual problem efficiently and in a distributed fashion. 

The remainder of the paper presents various results, such as upper and lower bounds for problems from the QAP and TSP libraries. We show that the proposed relaxation significantly reduces the running time compared to alternative SDP relaxations of the same complexity, while still producing strong lower and upper bounds. The assignment problem from Nuclear Magnetic Resonance Spectroscopy (NMR) is also formulated as a QAP problem, and results on benchmark synthetic datasets are presented which suggest that the new relaxation is a promising tool to tackle the problem, comparing favorably to state-of-the-art algorithms.

\subsection{Notation} Capitalized Roman letters, such as $A$, represent matrices, while their lower case equivalents stand for the corresponding column-wise vectorization, i.e. $a:=\vect{A}$. The symbol $\otimes$ is used to denote the Kronecker product. $\Pi_\mcal K$ represents a projection into the convex space $\mcal K$. $\Perm$ denotes the set of permutation matrices of dimension $n$ while $\text{DS}(n)$ denotes the set of doubly stochastic matrices. $I_n$ is the identity matrix of dimension $n$, $J_n$ is the all ones matrix of dimension $n$ and $\mbf 1_n$ is the all ones column vector of length $n$. For a matrix $Q$, the notation $Q_{ij}$ denotes the $(i,j)$-th block of the matrix (whose size should be clear from context), while $Q(i,j)$ denotes the $(i,j)$-th entry. Column $i$ of matrix $P$ is denoted by $p_i$, and $v(i)$ is used to indicate the $i$-th entry of vector $v$. Finally, we use $\delta_{ij}$ to denote the Kronecker-delta.
\section{Related Work}\label{sec:related}
Making use of the cyclic properties of the trace and of the vectorization identity $\vect{AYB}=(B^T\otimes A)\vect{Y}$, one can rewrite the QAP objective as follows:
\begin{align*}
\trace{PBP^TA^T} & = \trace{P^TAPB^T} \\
& = \vect{P}^T\vect{APB^T}\\
& = \vect{P}^T(B\otimes A)\vect{P}\\
& = \trace{(B\otimes A)\vect{P}\vect{P}^T}.
\end{align*}

The problem can therefore be reformulated as

\begin{problem}[Quadratic Assignment Problem]\label{prob:1}
	\begin{align*}
	\max_{Q,P} & \qquad \trace{(B\otimes A)Q}\\
	\st  & \qquad P\in\Perm\\
	& \qquad Q=\vect{P}\vect{P}^T.
	\end{align*}	
\end{problem}

Note that the constraints on $P$ and $Q$ are both nonconvex. One can relax $P$ to the set of doubly stochastic matrices. The nonconvex constraint on $Q$ can be replaced by $Q-\vect{P}\vect{P}^T\succeq 0$, which, by the Schur complement, is equivalent to:
\begin{equation*}
\left[\begin{array}{cc}
Q & \vect{P} \\
\vect{P}^T & 1
\end{array}\right]\succeq 0.
\end{equation*}

Enforcing additional linear constraints on $Q$ arising from the fact that each block $Q_{ij}$ of $Q$ is the outer product of two columns of a permutation matrix, one arrives at the convex relaxation proposed in \cite{Wolkowicz}:

\begin{problem}[SDP relaxation by Zhao \etal]
	\begin{align}
	\max_{Q,P} & \qquad \trace{(B\otimes A)Q}\\
	\st  & \qquad \left[\begin{array}{cc}
	Q & \vect{P} \\
	\vect{P}^T & 1
	\end{array}\right]\succeq 0,\\
	& \qquad \sum_i Q_{ii} = I_n,\ i=1,\ldots,n, \label{con:identity}\\
	& \qquad \trace{Q_{ij}} = 0,\ i\neq j,\ i,j=1,\ldots,n, \label{con:orth1}\\
	& \qquad \trace{Q_{ij}J_n} = 1,\ i,j=1,\ldots,n, \label{con:orth2}\\
	& \qquad Q_{ii}(j,j) =P(j,i),\ i=1,\ldots,n, \label{con:squared}\\
	& \qquad P \in \text{DS}(n),\  i=1,\ldots,n, \\
	& \qquad Q_{ij} \geq 0,\ i,j=1,\ldots,n.
	\end{align}
\end{problem}
Constraint \ref{con:identity} arises from the fact that each diagonal block of the un-relaxed $Q$ is the outer product of one of the columns of the permutation matrix with itself, such that it must have a single 1 on its diagonal. Constraints \ref{con:orth1} and \ref{con:orth2} arise from the orthogonality of the columns of a permutation matrix. Constraint \ref{con:squared} follows from the fact that the diagonal of $Q$ corresponds to the squared terms of the permutation, which are either $0$ or $1$, such that $Q_{ii}(j,j)=P(j,i)^2=P(j,i)$.

Although this relaxation is remarkably strong, achieving optimality in many problem instances in the QAP library, it is challenging to solve in practice. Problems of size $n>15$ are intractable on a regular computer using interior point methods, and even first-order methods are extremely slow for problems of size $n>50$.

\vspace{-3mm}
\section{Edge-based SDP relaxation and its generalization to clique-based SDP}\label{sec:csdp}

In many interesting applications, $B$ is a sparse matrix with $\mcal{O}(n)$ nonzero entries. This is the case for the traveling salesman problem (TSP) and longest-path problem. Denoting the set of edges in $\mcal{G}_B$ by $E(\mcal{G}_B)$, the QAP cost can be decomposed as
\begin{equation}
\trace{(B\otimes A)Q} = \sum_{i=1}^n \sum_{j=1}^n \trace{ B(i,j) A^T Q_{ij}} = \sum_{(i,j)\in E(\mathcal{G}_B)} \trace{ B(i,j) A^T Q_{ij}}.
\end{equation}
where $Q_{ij}$ is the $(i,j)$-th block of Q. The edge-SDP (E-SDP for short) relaxation we propose leverages the fact that in a graph where the adjacency matrix $B$ is sparse, the majority of the terms in $Q$ do not contribute to the objective function. Then a problem size-reduction is achieved by retaining only the blocks of $Q$ which are featured in the objective, leading to the following relaxation
\begin{problem}[E-SDP relaxation]
	Given a connected graph $B$ with $\mcal{O}(n)$ edges and an arbitrary graph $A$, solve
	\begin{align*}
	\max_{\{Q_{ij}\},P} & \qquad \sum_{(i,j)\in E(\mathcal{G}_B)} \trace{ B(i,j) A^T Q_{ij}} \\
	\st & \qquad \left[\begin{array}{ccc}
	Q_{ii} & Q_{ij} & p_i \\
	Q_{ij}^T & Q_{jj} & p_j \\
	p_i^T & p_j^T & 1
	\end{array}\right]\succeq 0,\ \forall i,j, \\
	& \qquad \trace{Q_{ij}} = 0,\ i\neq j,\ \\
	& \qquad \trace{Q_{ij}(J_n-I_n)} = 1,\ i\neq j,\ \\
	& \qquad \trace{Q_{ii}} = 1, \forall i,\ \\
	& \qquad \trace{Q_{ii}J_n} = 0, \forall i, \ \\
	& \qquad Q_{ii}(j,j) =P(j,i),\ \forall i,j, \\
	& \qquad P \in \text{DS}(n), \\
	& \qquad Q_{ij} \geq 0,\ \forall i,j,
	\end{align*}
\end{problem}
where we remind the reader that $p_i$ denotes the $i$-th column of matrix $P$. 

\paragraph{Semidefiniteness of Q:} For a large graph $B$ with $\mcal{O}(n)$ edges, we see that this relaxation has on the order of $\mcal{O}(n^3)$ variables. To achieve this reduction in problem size we sacrifice positive semidefiniteness of $Q$, and instead enforce only positive semidefiniteness of submatrices of $Q$, so this is a strictly weaker relaxation.

\subsection{Generalization to Clique-SDP (C-SDP)} \label{ssec:csdp}

In order to strengthen the relaxation, one can consider cliques of arbitrary size in $B$, rather than edges. Let $r$ denote such a clique in $B$ with nodes $V_r=\{V_r(1),\ldots,V_r({\vert V_r \vert})\}$. To simplify notation down the line, and because it is important in writing the ADMM formulation in section \ref{sec:admm}, we introduce the variable $X_r$ defined as follows:
\begin{equation}\label{eq:clique_var}
X_r = \left[\begin{array}{ccccc}
Q_{V_r(1)V_r(1)} & Q_{V_r(1)V_r(2)} & \cdots & Q_{V_r(1)V_r({\vert V_r \vert})} & p_{V_r(1)} \\
Q_{V_r(2)V_r(1)} & Q_{V_r(2)V_r(2)} & \cdots & Q_{V_r(2)V_r({\vert V_r \vert})} & p_{V_r(2)} \\
\vdots     & \vdots     & \ddots & \vdots     & \vdots \\
Q_{V_r({\vert V_r \vert}) V_r(1)} & Q_{V_r({\vert V_r \vert}) V_r(2)} & \cdots & Q_{V_r({\vert V_r \vert})V_r({\vert V_r \vert})} & p_{V_r({\vert V_r \vert})} \\
p_{V_r(1)}^T  & p_{V_r(2)}^T  & \cdots & p_{V_r({\vert V_r \vert})}^T  & 1
\end{array}\right]\succeq 0.
\end{equation}

With an appropriate choice of cost matrices $C_r$, and a variable $X_r$ for each clique in $B$, the QAP objective can be rewritten in terms of these variables as follows:
\begin{equation}
\trace{(B\otimes A)Q} = \sum_r\trace{C_rX_r}.
\end{equation}

The constraints for the edge-based case extend straightforwardly to the case of arbitrary cliques. An additional constraint arises from the fact that cliques share nodes, such that equality constraints need to be enforced between submatrices of variables $X_r$ corresponding to different cliques. This challenge is illustrated with an example in section $\ref{sec:path illustration}$ below, but we will see in section \ref{sec:admm} that formulating the problem in this manner greatly facilitates the development of an efficient and distributed ADMM scheme to solve the problem.

\paragraph{Computational issues:} In practice, it is computationally infeasible to consider all cliques in the graph. Instead, it is often worthwhile to consider only cliques of smaller size. Further, while the ADMM formulation can be used with mixed variable sizes, it is convenient to consider cliques of fixed size. To enforce this, one can define a maximal clique size and merge smaller cliques as needed to form variables of the required size. This same technique can be used to form variables of a fixed size on a dense subgraph.

\subsection{Illustrative example using the path graph}\label{sec:path illustration}
Let $B$ be the adjacency matrix for the path graph, given by
\begin{equation*}
B = \begin{bmatrix}
0      & 1      & \cdots & \cdots & 0 \\
\vdots & & \ddots &  & \vdots \\
\vdots &        &  &  & 1 \\
0      & \cdots & \cdots & \cdots & 0
\end{bmatrix}
\end{equation*}

In this case, the cliques are the $n-1$ edges of $B$. As a result, we will have $n-1$ variables of the form
\begin{equation*}
X_{i (i+1)} = \left[\begin{array}{ccc}
Q_{ii} & Q_{i(i+1)} & p_i \\
Q_{(i+1)i} & Q_{(i+1)(i+1)} & p_{i+1} \\
p_i^T & p_{i+1}^T & 1
\end{array}\right].
\end{equation*}

If we consider an example with 5 nodes and look at the matrix $Q=\vect{P}\vect{P}^T$, we see that the E-SDP variables include the $2n\times 2n$ blocks along the diagonal, as illustrated in Figure \ref{fig:e-sdp-2}.

\begin{figure}
	\centering
	\begin{subfigure}[b]{0.45\textwidth}
		\begin{center}
			\begin{tikzpicture}
			\matrix [matrix of math nodes,left delimiter=(,right delimiter=), ampersand replacement=\&] (m)
			{
				Q_{11} \& Q_{12} \& Q_{13} \& Q_{14} \& Q_{15} \\
				Q_{21} \& Q_{22} \& Q_{23} \& Q_{24} \& Q_{25} \\
				Q_{31} \& Q_{32} \& Q_{33} \& Q_{34} \& Q_{35} \\
				Q_{41} \& Q_{42} \& Q_{43} \& Q_{44} \& Q_{45} \\
				Q_{51} \& Q_{52} \& Q_{53} \& Q_{54} \& Q_{55} \\
			};
			\draw[color=red,opacity=1] (m-1-1.north west) -- (m-1-2.north east) -- (m-2-2.south east) -- (m-2-1.south west) -- (m-1-1.north west);
			\draw[color=red,opacity=1] (m-2-2.north west) -- (m-2-3.north east) -- (m-3-3.south east) -- (m-3-2.south west) -- (m-2-2.north west);
			\draw[color=red,opacity=1] (m-3-3.north west) -- (m-3-4.north east) -- (m-4-4.south east) -- (m-4-3.south west) -- (m-3-3.north west);
			\draw[color=red,opacity=1] (m-4-4.north west) -- (m-4-5.north east) -- (m-5-5.south east) -- (m-5-4.south west) -- (m-4-4.north west);
			\begin{scope}[on background layer]
			\fill[color=blue!30,opacity=1] (m-2-2.north west) rectangle (m-2-2.south east);
			\fill[color=blue!30,opacity=1] (m-3-3.north west) rectangle (m-3-3.south east);
			\fill[color=blue!30,opacity=1] (m-4-4.north west) rectangle (m-4-4.south east);
			\end{scope}
			\node [color=blue, above, opacity=1] at (m-1-1.north east) {$X_{12}$(1:2$n$,1:2$n$)};
			\node [color=blue, below, opacity=1] at (m-5-4.south east) {$X_{45}$(1:2$n$,1:2$n$)};
			\end{tikzpicture}
		\end{center}
		\caption{E-SDP with cliques of size 2.}
		\label{fig:e-sdp-2}
	\end{subfigure}
	\begin{subfigure}[b]{0.45\textwidth}
		\begin{center}
			\begin{tikzpicture}
			\matrix [matrix of math nodes,left delimiter=(,right delimiter=), ampersand replacement=\&] (m)
			{
				Q_{11} \& Q_{12} \& Q_{13} \& Q_{14} \& Q_{15} \\
				Q_{21} \& Q_{22} \& Q_{23} \& Q_{24} \& Q_{25} \\
				Q_{31} \& Q_{32} \& Q_{33} \& Q_{34} \& Q_{35} \\
				Q_{41} \& Q_{42} \& Q_{43} \& Q_{44} \& Q_{45} \\
				Q_{51} \& Q_{52} \& Q_{53} \& Q_{54} \& Q_{55} \\
			};
			\draw[color=red,opacity=1] (m-1-1.north west) -- (m-1-3.north east) -- (m-3-3.south east) -- (m-3-1.south west) -- (m-1-1.north west);
			\draw[color=red,opacity=1] (m-2-2.north west) -- (m-2-4.north east) -- (m-4-4.south east) -- (m-4-2.south west) -- (m-2-2.north west);
			\draw[color=red,opacity=1] (m-3-3.north west) -- (m-3-5.north east) -- (m-5-5.south east) -- (m-5-3.south west) -- (m-3-3.north west);
			\begin{scope}[on background layer]
			\fill[color=blue!30,opacity=1] (m-2-2.north west) rectangle (m-2-2.south east);
			\fill[color=blue!30,opacity=1] (m-3-3.north west) rectangle (m-3-3.south east);
			\fill[color=blue!30,opacity=1] (m-3-3.north west) rectangle (m-4-4.south east);
			\fill[color=orange!30,opacity=1] (m-2-3.north west) rectangle (m-2-3.south east);
			\fill[color=orange!30,opacity=1] (m-3-2.north west) rectangle (m-3-2.south east);
			\fill[color=orange!30,opacity=1] (m-3-4.north west) rectangle (m-3-4.south east);
			\fill[color=orange!30,opacity=1] (m-4-3.north west) rectangle (m-4-3.south east);
			\end{scope}
			\node [color=blue, above, opacity=1] at (m-1-1.north east) {$X_{123}$(1:3$n$,1:3$n$)};
			\node [color=blue, below, opacity=1] at (m-5-4.south east) {$X_{345}$(1:3$n$,1:3$n$)};
			\end{tikzpicture}
		\end{center}
		\caption{C-SDP with cliques of size 3.} \label{fig:c-sdp-3}
	\end{subfigure}
	\caption{C-SDP variables in the path graph problem with 5 nodes.}
\end{figure}


Several of the diagonal blocks of $Q$, highlighted in blue, overlap between adjacent variables, and thus it is necessary to enforce these equality constraints. Fortunately, the diagonal blocks of $Q$ are diagonal themselves, such that $n$ equality constraints are sufficient to enforce equality between two blocks.

We draw attention to the fact that when using cliques of size greater than 2 the C-SDP variables $X_r$ will overlap in off-diagonal blocks (as illustrated in Figure \ref{fig:c-sdp-3}). These are problematic to handle computationally (as they are generally dense), but, in practice, we have observed that enforcing equalities only between diagonal blocks sacrifices little in terms of performance, while greatly reducing the computational burden, as we shall see from the ADMM scheme in section \ref{sec:admm} below.
\section{Alternating Direction Method of Multipliers}\label{sec:admm}
In this section, we devise an ADMM that solves the E-SDP and C-SDP relaxations in a distributed manner. We present the updates in each ADMM iteration for E-SDP, and this extends to C-SDP in a straightforward manner.

\subsection{Rewriting constraints in E-SDP relaxation}
We saw in section \ref{ssec:csdp} how the E-SDP objective can be written in terms of the variables $X_{ij}$ and cost matrices $C_{ij}$. 

\paragraph{Note on convention:} Since the variables $X_{ij}$ are symmetric, it is equivalent to consider $X_{ij}$ or $X_{ji}$. Therefore, we define graph $\mathcal{G}_{\tilde{B}}$ with adjacency matrix, $\tilde B$ defined as
\begin{equation}
\tilde{B}=\text{triu}(B+B^T)
\end{equation}
where $\text{triu}(M)$ extracts the upper triangular portion of matrix $M$. Thus, in all that follows, the variables $X_{ij}$ will be defined according to the edges specified by $\tilde{B}$, such that $i<j$.

It remains to rewrite the constraints in terms of these variables. We remind the reader that the variables under consideration are
\begin{equation}
X_{ij} = \left[\begin{array}{ccc}
Q_{ii} & Q_{ij} & p_i \\
Q_{ij}^T & Q_{jj} & p_j \\
p_i^T & p_j^T & 1
\end{array}\right],\quad (i,j)\in E(\mathcal{G}_{\tilde{B}})
\end{equation}
where all $X_{ij}$'s are non-negative and PSD. Going forward, all the constraints will be rewritten in terms of the variables $X_{ij}$ and the variables $Q_{ij}$ and $p_i$ will no longer be used. 

The first set of constraints on $X_{ij}$ follows directly from the constraints on $Q$ and $P$, which are the following:
\begin{align}
& \trace{Q_{ii}} = \trace{Q_{jj}} = 1, \\
& \trace{Q_{ii}J_n} = \trace{Q_{jj}J_n} = 0, \\
& \trace{Q_{ij}} = 0, \\
& \trace{Q_{ij}J_n} = 1, \\
& \text{diag}(Q_{ii}) = p_i, \label{eq:last_row_col1} \\
& \text{diag}(Q_{jj}) = p_j. \label{eq:last_row_col2}
\end{align}
Letting $x_{ij}:=\vect{X_{ij}}$, all the constraints above can be written in the form:
\begin{equation}\label{eq:lin_constraints}
\mcal{A}x_{ij} = b_E
\end{equation}
where $\mcal{A}\in \mathbb{R}^{m_E\times (2n+1)^2}$, $b_E\in \mathbb{R}^{m_E}$, and $m_E$ is the number of equality constraints. 

Note that the $X_{ij}$'s are not independent of each other. Firstly, for the edges that are incident on the same node, the associated variables $X_{ij}$'s share a common $n\times n$ block on the diagonal. This is illustrated in the example of a path graph in section \ref{sec:path illustration}. Therefore, equality constraints between the overlapping diagonal blocks of $X_{ij}$'s have to be enforced. Since $\trace{Q_{ii}(J_n-I_n)} = \trace{Q_{jj}(J_n-I_n)} = 0$ and $Q_{ii},Q_{jj}\geq 0$, the off-diagonal terms of $Q_{ii}$ and $Q_{jj}$ are zeros and it suffices to enforce equality of the diagonals. Further, since $p_i$ and $p_j$ equal the diagonals of $Q_{ii}$ and $Q_{jj}$, one can enforce consistency of the overlapping blocks by looking at the last row and column of each $X_{ij}$ instead. Consider the sampling matrices $\mcal{B}_1$ and $\mcal{B}_2$, which sample $p_i$ and $p_j$ from the vector $x_{ij}$ above. If $(i,j),(k,i)\in E(\mathcal{G}_{\tilde{B}})$, then a consistency relationship of the form
\begin{equation}
\label{consistency relation}
\mathcal{B}_1 x_{ij} = \mathcal{B}_2 x_{ki}
\end{equation}
must hold.

Adding the conic constraints for positivity and positive semi-definiteness, the E-SDP relaxation can be reformulated as:
\begin{problem} \label{prob:4}
	\begin{align*}
	\max_{\{x_{ij}\}} & \qquad \sum_{(i,j)\in E(\mathcal{G}_{\tilde{B}})} c_{ij}^Tx_{ij} \\
	\st & \qquad \mcal A x_{ij} = b_e, \\
	& \qquad \mcal B_1 x_{ij} = \mcal B_2 x_{ki}\,  \forall (i,j), (k,i)\in E(\mcal{G}_{\tilde{B}}),\ i=1,\ldots,n, \\
	& \qquad \mcal{B}_2x_{kn}+\sum_{i=1}^{n-1} \mcal{B}_1x_{ij_i}=\mbf 1, (i,j_i), (k,n) \in E(\mcal{G}_{\tilde{B}}) \\
	& \qquad x_{ij} \succeq 0, \\
	& \qquad \mcal D x_{ij} \geq 0.
	\end{align*}
\end{problem}
Here, with a slight abuse of notation, we have used $x_{ij} \succeq 0$ to denote $X_{ij}\succeq 0$. We have also used $c_{ij}\equiv \vect{C_{ij}}$. The third constraint amounts to stating that the sum of the diagonal blocks of $Q$ equal the identity. The matrix $\mcal{D}$ is of size $4n^2\times (2n+1)^2$ and it samples all elements of $x_{ij}$ except for those corresponding to the last row and column of $X_{ij}$. The reason for this sampling is two-fold:
\begin{enumerate}
	\item Sampling the last row and column is unnecessary, since these entries are implicitly defined by the linear constraints covered in equations \ref{eq:last_row_col1} and \ref{eq:last_row_col2};
	\item Using a sampling operator of this form ensures mutual orthogonality between $\mcal{D}$, $\mcal{B}_1$ and $\mcal{B}_2$, 
	\begin{equation}\label{eq:mutual_orth}
	\mcal{B}_1\mcal{B}_2^T=\mathbf{0}_{n\times n},\qquad \mcal{B}_1\mcal{D}^T=\mathbf{0}_{n\times 4n^2}, \qquad\mcal{B}_2\mcal{D}^T=\mathbf{0}_{n\times 4n^2},
	\end{equation} which shall prove crucial in obtaining fast ADMM updates involving a least-squares problem with a block-diagonalized Hessian.
\end{enumerate}

 In order to derive a fast ADMM routine to solve Problem \ref{prob:4}, slack variables are introduced. Let $N_i$ be the one-hop neighborhood of node $i$ on graph $\mcal{G}_{\tilde{B}}$. Then the consistency relation in equation \ref{consistency relation} can be enforced by introducing slack variables $p_i$, such that
\begin{align}
& \mcal{B}_1x_{ij} = p_i,\ \forall j\in N_i\\
& \mcal{B}_2x_{ij} = p_j,\ \forall j\in N_i.
\end{align}

As a result, our problem can finally be written in the form
\begin{problem}   \label{prob:5}
	\begin{align*}
	\max_{\{x_{ij}\},\{p_i\}} & \qquad \sum_{(i,j)\in E(\mathcal{G}_{\tilde{B}})} c_{ij}^Tx_{ij} \\
	\st & \qquad y_{ij} :\mcal A x_{ij} = b_e, \\
	& \qquad w_{ij}^{(1)} :\mcal B_1 x_{ij} = p_i\  \forall j\in N_i,\ i=1,\ldots,n, \\
	& \qquad w_{ij}^{(2)} :\mcal B_2 x_{ij} = p_j\  \forall j\in N_i,\ i=1,\ldots,n, \\
	& \qquad t : \sum_i p_i = \bf{1}, \\
	& \qquad s_{ij}\succeq 0 : x_{ij} \succeq 0, \\
	& \qquad z_{ij}\geq 0 : \mcal D x_{ij} \geq 0
	\end{align*}
\end{problem}
where the variable in front of each colon is the dual variable tied to the corresponding constraint. The constraint $\sum_{i=1}^n p_{i} = \mathbf{1}_n$ couples the $X_{ij}$ from different blocks together. 

\paragraph{Generalization to C-SDP}: One can generalize the presentation above to the clique-based SDP relaxation in a straightforward way. The one important difference is that for general sets of nodes of size greater than two, the corresponding $X_r$'s as defined by equation \ref{eq:clique_var} might overlap in non-diagonal blocks (if two or more nodes are shared by two cliques).

Figure \ref{fig:c-sdp-3} highlights the fact that one must enforce equalities between $X_{ijk}$ and $X_{ijl}$ not only for $Q_{ii}$ and $Q_{jj}$, but also for $Q_{ij}$. However, we have observed that enforcing only the equalities on the diagonal blocks produces solutions which are nearly as good with a much decreased computational cost, so we adopt this solution for the remainder of the paper

\subsection{Dual problem and the ADMM updates}
We now turn to the dual problem of the E-SDP relaxation presented in the form of problem \ref{prob:5}. In this section we show that with a proper grouping of the dual variables the ADMM updates for solving the dual problem can be computed in a distributed manner.

The dual of problem \ref{prob:5} is the following:
\begin{problem}\label{prob:6}
	\begin{align*}
	\min_{y_{ij},w_{ij}^{(k)},t,s_{ij},z_{ij}} &\qquad\sum_{(i,j)\in E(\mcal{G}_{\tilde{B}})}b_e^Ty_{ij}-\mbf{1}^Tt \\
	\st & \qquad s_{ij}\succeq 0, (i,j)\in E(\mcal{G}_{\tilde{B}}) \\
	& \qquad z_{ij}\geq 0, (i,j)\in E(\mcal{G}_{\tilde{B}})  \\
	& \qquad x_{ij}:-c_{ij}+s_{ij}+\mcal{D}^Tz_{ij}+\mcal{A}^Ty_{ij}+\\
	& \qquad \qquad \qquad \qquad \mcal{B}_1^Tw_{ij}^{(1)}+\mcal{B}_2^Tw_{ij}^{(2)} = 0, (i,j)\in E(\mcal{G}_{\tilde{B}}) \\
	& \qquad g_i: t-\sum_{j\in N_i} w_{ij}^{(1)}-\sum_{j:i \in N_j} w_{ji}^{(2)} = 0, i=1,\ldots, n.
	\end{align*}
\end{problem}

Using $\delta_{\mcal{K}}(x)$ to denote a function that takes the value $+\infty$ for $x\notin \mcal{K}$ and $0$ otherwise, the augmented Lagrangian is  
\begin{multline}
\mcal{L} =  \sum_{(i,j)\in E(\mcal{G}_{\tilde{B}})} \left( \delta_{\mcal{S}_n^+}(s_{ij})+\delta_{\mcal{K}_p}(z_{ij})-b_e^Ty_{ij} \right)- \mbf 1^Tt +\\
\frac{\rho}{2}\sum_{(i,j)\in E(\mcal{G}_{\tilde{B}})} \bnorm{-c_{ij}+s_{ij}+\mcal D^T z_{ij}+\mcal{A}^Ty_{ij}+ \mcal B_1^T w_{ij}^{(1)} + \mcal{B}_2^Tw_{ij}^{(2)}+\frac{x_{ij}}{\rho}}_2^2+\\
\frac{\rho}{2}\sum_{i=1}^n \bnorm{t- \sum_{j\in N_i}   w_{ij}^{(1)} - \sum_{j:i \in N_j}  w_{ji}^{(2)}+\frac{g_i}{\rho}}^2_2
\end{multline}
where $x_{ij}$ and $g_i$ are now the dual variables of the dual problem, and $\rho$ is some constant greater than 0.

In \cite{Toh}, a convergent ADMM is proposed to solve optimization problems with a 3-block structure where one of the blocks only involves linear operators. In our problem, we let the three blocks be defined by the groups of variables $(s_{ij},t)$, $(y_{ij})$ and $(z_{ij},w_{ij}^{(k)})$. Then the algorithm proceeds as follows:

\begin{algorithm}
	\caption{Conic-ADMM3c \protect\cite{Toh}}
	\begin{algorithmic}
		\REQUIRE $\rho>0$	and $\tau=1$
		\FOR{$l=1,\ldots,\text{MAXIT}$}
		\item $(s_{ij},t)^{l+1}\leftarrow\argmin_{s_{ij},t}\mcal L(s_{ij},t,y^l_{ij},z^l_{ij},w_{ij}^{(k),l}; x_{ij}^l,g_i^l;\rho)$
		\item $(y_{ij})^{l+1/2}\leftarrow\argmin_{y_{ij}}\mcal L(s_{ij}^{l+1},t^{l+1},y_{ij},z^l_{ij},w_{ij}^{(k),l}; x_{ij}^l,g_i^l;\rho)$
		\item $(z_{ij},w_{ij}^{(k)})^{l+1}\leftarrow\argmin_{z_{ij},w_{ij}^{(k)}}\mcal L(s_{ij}^{l+1},t^{l+1},y^{l+1/2}_{ij},z_{ij},w_{ij}^{(k)}; x_{ij}^l,g_i^l;\rho)$
		\item $(y_{ij})^{l+1}\leftarrow\argmin_{y_{ij}}\mcal L(s_{ij}^{l+1},t^{l+1},y_{ij},z^{l+1}_{ij},w_{ij}^{(k),l+1}; x_{ij}^l,g_i^l;\rho)$
		\item $x_{ij}^{l+1}\leftarrow x_{ij}^l+\tau\argmin_{x_{ij}} \mcal{L}(s_{ij}^{l+1},t^{l+1},y_{ij}^{l+1},z^{l+1}_{ij},w_{ij}^{(k),l+1}; x_{ij},g_i^l;\rho)$
		\item $g_i^{l+1}\leftarrow g_i+\tau\argmin_{g_{i}} \mcal{L}(s_{ij}^{l+1},t^{l+1},y_{ij}^{l+1},z^{l+1}_{ij},w_{ij}^{(k),l+1}; x_{ij}^{l+1},g_i;\rho)$
		\ENDFOR
	\end{algorithmic}
\end{algorithm}

In the remainder of this section, we will derive each of the updates in turn and illustrate how this choice of variable groupings allows for easy parallelization.

\paragraph{Update for $(s_{ij},t)$:} The updates for $s_{ij}$ and for $t$ are independent. The update for $t$ is given by the solution to a least-squares problem:
\begin{equation}
\argmin_{t} \mathcal{L} = \frac{1}{n}\left(\sum_{i=1}^n\sum_{j\in N_i}   w_{ij}^{(1)} + \sum_{i=1}^n\sum_{j:i \in N_j}  w_{ji}^{(2)}-\frac{1}{n\rho}\sum_ig_i\right)+\frac{1}{n\rho}\mbf 1.
\end{equation}
The new $s_{ij}$ is obtained from
\begin{equation}
\argmin_{s_{ij}} \mathcal{L}= 
\Pi_{\mcal S_n^+}(c_{ij}-\mcal D^Tz_{ij}-\mcal A^Ty_{ij}-\mcal{B}_1^Tw_{ij}^{(1)}-\mcal{B}_2^Tw_{ij}^{(2)}-\rho^{-1}x_{ij}),
\end{equation}
where $\Pi_{\mcal{S}_n^+}$ is a projection to the positive semidefinite cone.

\paragraph{Update for $(y_{ij})$:} The update for $y_{ij}$ is the solution to a least-squares problem, given by
\begin{equation}
\argmin_{y_{ij}} \mathcal{L}= 
(\mcal A\mcal A^T)^{-1}(\mcal A(c_{ij}-s_{ij}-\mcal D^Tz_{ij}  -\mcal B_1^Tw_{ij}^{(1)}-\mcal B_2^Tw_{ij}^{(2)}-\rho^{-1}x_{ij})+\rho^{-1}b_e).
\end{equation}

By construction, $\mcal{A}$ is the matrix that encodes the linear constraints. Note that $\mcal{A}$ is of size $m_E\times \mcal{O}(n^2)$ and has linearly independent rows. Since $m_E$ is of order $\mcal{O}(n)$,  $\mcal{A}\mcal{A}^T$ is a full-rank matrix of dimension $\mcal{O}(n)$.

\paragraph{Update for $(z_{ij},w_{ij}^{(k)})$:} The updates for $z_{ij}$ and for the $w_{ij}^{(k)}$'s decouple due to the fact that the sampling matrices $\mcal D$ and $\mcal B_1$ and $\mcal B_2$ have mutually orthogonal rows, as they sample different entries of $X_{ij}$. To see this, we write the relevant minimization problem as follows
\begin{multline*}
\min_{w_{ij}^{(k)},z_{ij}} \sum_{(i,j)\in E(\mcal{G}_{\tilde{B}})} \delta_{\mcal{K}_p}(z_{ij})
\\+\frac{\rho}{2}\sum_{(i,j)\in E(\mcal{G}_{\tilde{B}})} \bnorm{-c_{ij}+s_{ij}+\mcal D^T z_{ij}+\mcal{A}^Ty_{ij}+ \mcal B_1^T w_{ij}^{(1)} + \mcal{B}_2^Tw_{ij}^{(2)}+\frac{x_{ij}}{\rho}}_2^2\\
+\frac{\rho}{2}\sum_{i=1}^n \bnorm{t- \sum_{j\in N_i}   w_{ij}^{(1)} - \sum_{j:i \in N_j}  w_{ji}^{(2)}+\frac{g_i}{\rho}}^2_2.
\end{multline*}
Recalling that $N_i$ is the set of one-hop neighbors of node $i$, define
\begin{equation*}
\mcal K_{N_i}=\left[\begin{array}{ccccccccc}
\mcal B_1^T & & & \mcal B_2^T & & & \mcal D^T & & \\
& \ddots & & & \ddots & & & \ddots & \\
& & \mcal B_1^T & &  & \mcal B_2^T & & & \mcal D^T \\
I & \cdots & I & I & \cdots & I & &
\end{array}\right]
\end{equation*}
which has $|N_i|$ copies of $\mcal{B}_1^T$, $\mcal{B}_2^T$ and $\mcal{D}^T$, $2|N_i|$ copies of the identity, and is zero everywhere else. Then the problem becomes
\begin{multline} \label{eq:min_wz}
\min_{w_{ij}^{(k)},z_{ij}} \sum_{(i,j)\in E(\mcal{G}_{\tilde{B}})}\delta_{\mcal{K}_p}(z_{ij})  +\sum_{i=1}^n\frac{\rho}{2}\bnorm{\mcal{K}_{N_i}\text{vec}\left(V_{wz(i)}\right)+\left[\begin{array}{c}
	v_{iN_i(1)} \\ \vdots \\ v_{iN_i(|N_i|)} \\ -t-\rho^{-1}g_i
	\end{array}\right]}_2^2
\end{multline}
where
\begin{equation*}
V_{wz(i)}=\left[\begin{array}{ccccccccc}
w_{iN_i(1)}^{(1)} & \cdots & w_{iN_i(|N_i|)}^{(1)} & w_{iN_i(1)}^{(2)} & \cdots & w_{iN_i(|N(i)|)}^{(2)} & z_{iN_i(1)} & \cdots & z_{iN_i(|N_i|)}
\end{array}\right]
\end{equation*}
and
\begin{equation*}
v_{ij} = -c_{ij}+s_{ij}+\mcal{A}^Ty_{ij}+\frac{x_{ij}}{\rho}.
\end{equation*}
$\mcal{K}^T_{N_i}\mcal{K}_{N_i}$ has a block-diagonal structure, owing to the mutual orthogonality of the following blocks
\begin{equation*}
\left[\begin{array}{ccc}
\mcal B_1^T & & \\
& \ddots & \\
& & \mcal B_1^T \\
I & \cdots & I
\end{array}\right], \left[\begin{array}{ccc}
\mcal B_2^T & & \\
& \ddots & \\
& & \mcal B_2^T \\
I & \cdots & I
\end{array}\right], \left[\begin{array}{ccc}
\mcal D^T & & \\
& \ddots & \\
& & \mcal D^T \\
& & 
\end{array}\right],
\end{equation*}
as specified in \ref{eq:mutual_orth}.

The form of \ref{eq:min_wz} also shows that the problem can be solved independently for each neighborhood. The optimal value for $z_{ij}$ is then
\begin{equation}
z_{ij}=\Pi_{\mcal K_p}\left(\mcal D\left(c_{ij}-s_{ij}-\mcal A^Ty_{ij}-\frac{x_{ij}}{\rho}\right)\right),
\end{equation}
where $\Pi_{\mcal{K}_p}$ is a projection to the positive cone.

For $w_{ij}^{(k)}$, let
\begin{equation*}
\mcal B_{N_i}=\left[\begin{array}{cccccc}
\mcal B_1^T & & & \mcal B_2^T & &  \\
& \ddots & & & \ddots &  \\
& & \mcal B_1^T & &  & \mcal B_2^T\\
I & \cdots & I & I & \cdots & I
\end{array}\right].
\end{equation*}
Then
\begin{multline*}
\text{vec}\left(\left[\begin{array}{cccccc}
w_{iN_i(1)}^{(1)} & \cdots & w_{iN_i(|N_i|)}^{(1)} & w_{iN_i(1)}^{(2)} & \cdots & w_{iN_i(|N_i|)}^{(2)}
\end{array}\right]\right)=\\(\mcal B_{N_i}^T\mcal B_{N_i}^{})^{-1}B_{N_i}^T\left(\left[\begin{array}{c}
v_{iN_i(1)} \\ \vdots \\ v_{iN_i(|N_i|)} \\ -t-\rho^{-1}g_i
\end{array}\right]\right).
\end{multline*}
Note that the matrices $\mcal{B}_{N_i}^T\mcal{B}_{N_i}$ for $i=1,\ldots,n$ take the generic form
\begin{equation*}
\mcal{B}_{N_i}^T\mcal{B}_{N_i}=(\alpha-\beta) I_{|N_i|n} + \beta J_{|N_i|}\otimes I_n
\end{equation*}
where $\alpha$ and $\beta$ are constants. Therefore, their inverse is given by the $n|N_i|\times n|N_i|$ matrix
\begin{equation*}
(\mcal{B}_{N_i}^T\mcal{B}_{N_i})^{-1} = \frac{1}{\alpha-\beta}I_{|N_i|n}-\frac{\beta}{(\alpha-\beta)(\alpha-\beta+|N_i|\beta)} J_{|N_i|}\otimes I_n.
\end{equation*}
Let
\begin{equation}
\mcal{H}_i = \frac{1}{\alpha-\beta}I_{|N_i|}-\frac{\beta}{(\alpha-\beta)(\alpha-\beta+|N_i|\beta)} J_{|N_i|}
\end{equation}
which is a $|N_i|\times |N_i|$ matrix. Then, $(\mcal{B}_{N_i}^T\mcal{B}_{N_i})^{-1}v = v\mcal{H}_i$ and the update for $w_{ij}^{(k)}$'s is finally given by
\begin{equation}
\text{vec}\left(\left[\begin{array}{cccccc}
w_{iN_i(1)}^{(1)} & \cdots & w_{iN_i(|N_i|)}^{(1)} & w_{iN_i(1)}^{(2)} & \cdots & w_{iN_i(|N_i|)}^{(2)}
\end{array}\right]\right)=\mcal B_{N_i}^T\left(\left[\begin{array}{c}
v_{iN_i(1)} \\ \vdots \\ v_{iN_i(|N_i|)} \\ -t-\rho^{-1}g_i
\end{array}\right]\right)\mcal{H}_i.
\end{equation}

The updates for $x_{ij}$ and $g_i$, taken directly from \cite{Toh}, are the following:
\begin{equation}
x_{ij}^{k+1} = x_{ij}^k+\tau\rho\left(-c_{ij}+s_{ij}+\mcal D^Tz_{ij}+\mcal A^Ty_{ij}+\mcal B_1^Tw_{ij}^{(1)}+\mcal B_2^Tw_{ij}^{(2)}\right)
\end{equation}
and
\begin{equation}
g_{i}^{k+1}=g_{i}^k+\tau\rho\left(t- \sum_{j\in N_i}   w_{ij}^{(1)} - \sum_{j:i \in N_j}  w_{ji}^{(2)}\right).
\end{equation}

This concludes the ADMM formulation for Problem \ref{prob:6} using two-cliques (i.e. edges). The problem for larger cliques is very similar, with additional variables $w_{ij}^{(k)}$ for $k>2$.

The most costly update in the ADMM is a projection of a matrix of dimension $\mcal{O}(n)$ to the positive semidefinite cone. The updates can be parallelized (across the nodes, for e.g.), only requiring one gather operation per iteration for the $w_{ij}^{(k)}$ updates.

\subsection{Convergent ADMM vs. direct extension}
A convergent 3-block ADMM algorithm was not available until the paper by Sun \etal \cite{Toh}, yet it is common practice to use a direct extension of the 2-block ADMM algorithm, i.e. updating blocks in 1-2-3 order. In \cite{Chen}, the authors show that this straightforward extension is not necessarily convergent.

\begin{figure}
	\centering
	\begin{subfigure}[b]{0.45\textwidth}
		\includegraphics[width=\textwidth]{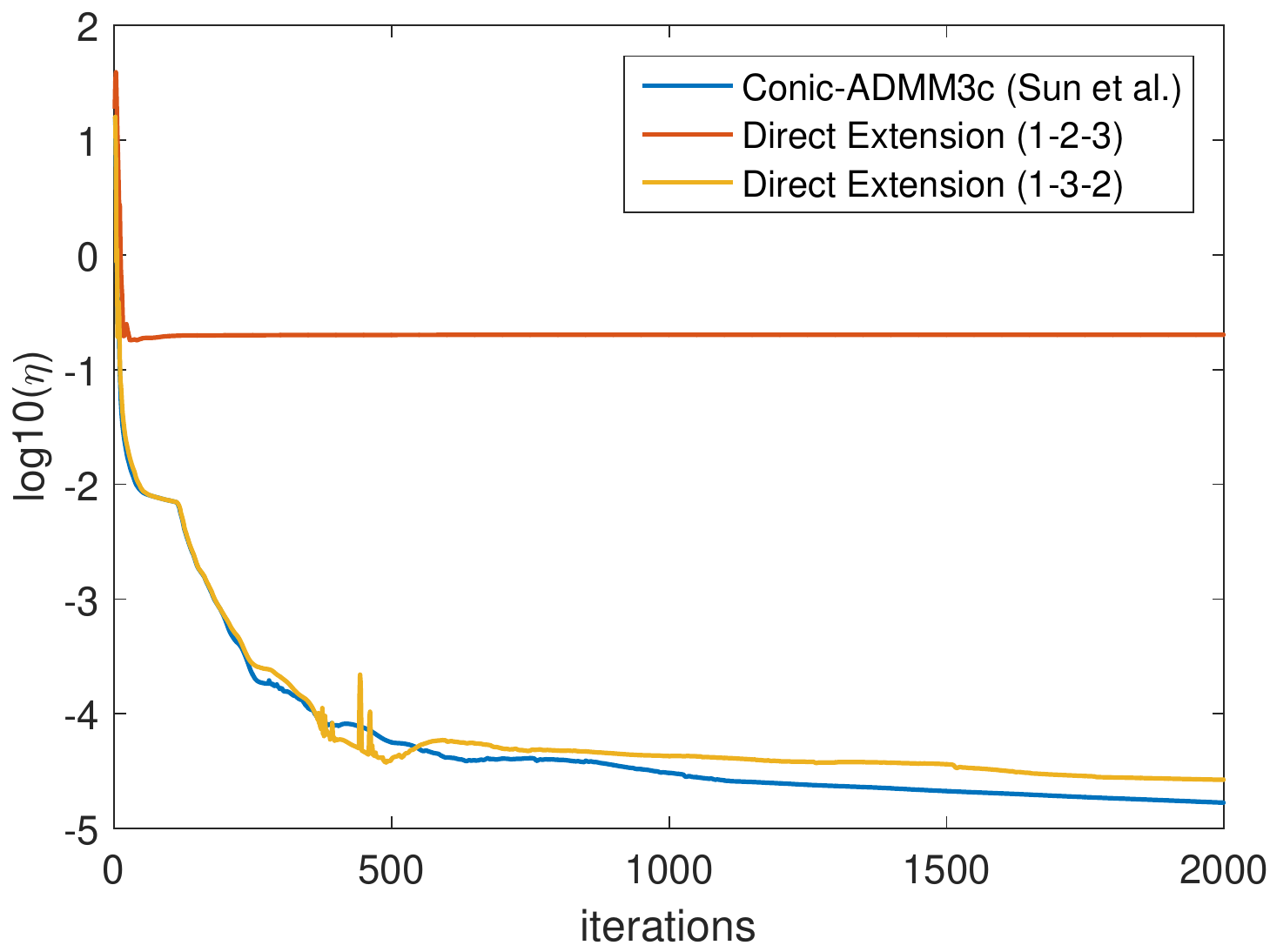}
		\caption{Problem gr21 (TSPLIB).}
		\label{fig:cvd_gr21}
	\end{subfigure}
	\begin{subfigure}[b]{0.45\textwidth}
		\includegraphics[width=\textwidth]{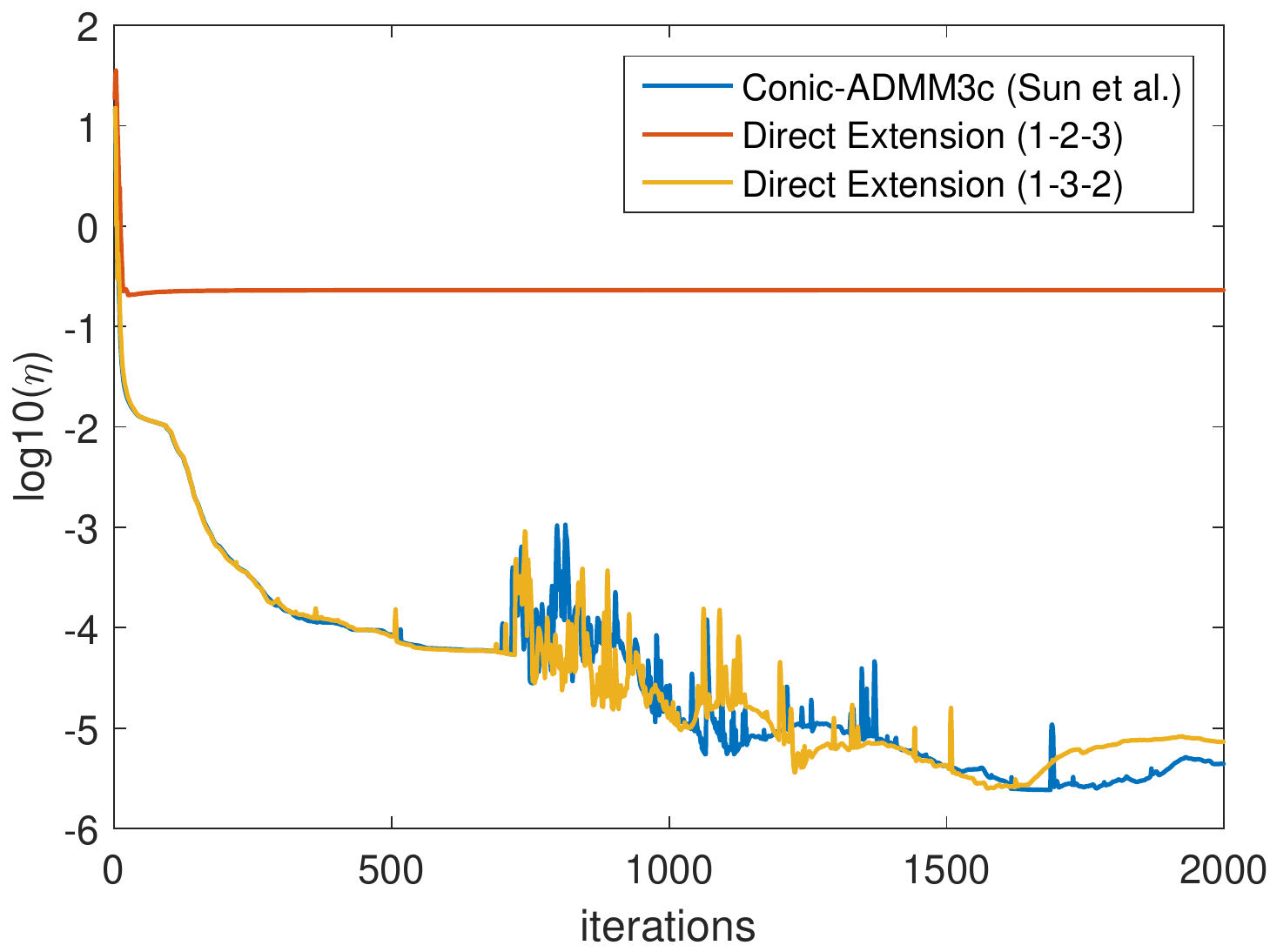}
		\caption{Problem chr20a (QAPLIB).}
		\label{fig:cvd_chr20a}
	\end{subfigure}
	\caption{Convergence criterion \ref{eq:conv} as a function of the number of iterations for problems from the TSPLIB and QAPLIB for both the Conic-ADMM3c algorithm of \protect\cite{Toh} (in blue) and the direct extension of 2-block ADMM to a multi-block setting (red and yellow). Updating blocks in order 1-2-3 fails to converge. Updating in order 1-3-2 yields a convergence curve similar to the one obtained by using the algorithm in \protect\cite{Toh}.}\label{fig:animals}
\end{figure}

Indeed, in this work, we have observed that a direct extension can fail to converge in a rather dramatic manner. As an example, Figures \ref{fig:cvd_gr21} and \ref{fig:cvd_chr20a} depict convergence curves for the gr21 problem (TSPLIB) and the chr20a problem (QAPLIB), respectively, ran to 2000 iterations for both Conic-ADMM3c and direct extension. When directly extending 2-block ADMM to 3-block ADMM with blocks $(s_{ij}, t)$, $(y_{ij})$, and $(z_{ij}, w_{ij}^{(k)})$, performing the updates in this order (1-2-3) fails to converge (although we observe convergence when updating in order 1-3-2).

Throughout our results, we make use of the convergence criterion $\eta$, defined analogously to the one in \cite{Toh}
\begin{equation}\label{eq:conv}
\eta = \max (\eta_P, \eta_D, \eta_{\mcal K}, \eta_{\mcal K*}, \eta_{\mcal P}, \eta_{\mcal P*}, \eta_{C1}, \eta_{C2})
\end{equation}
with 
\begin{align*}
\begin{array}{ll}
	\eta_P = \frac{\norm{\mcal A X-B_e}_F}{1+\sqrt{n}\norm{b_e}} & \qquad \eta_D = \frac{\norm{-C+\mcal A^TY+S+\mcal D^TZ+\mcal B_1^TW^{(1)}+B_2^TW^{(2)}}_F}{1+\sqrt{n}\norm{b_e}} \\
	\eta_{\mcal K}=\frac{\norm{\Pi_{\mcal S_+^n}(-X)}_F}{1+\norm{X}_F} &\qquad \eta_{\mcal K*} = \frac{\norm{\Pi_{\mcal S_+^n}(-S)}_F}{1+\norm{S}_F} \\
	\eta_{\mcal P} = \frac{\norm{X-\Pi_{\mcal K_p}(X)}_F}{1+\norm{X}_F} &\qquad \eta_{\mcal P*} = \frac{\norm{Z-\Pi_{\mcal K_p}(Z)}_F}{1+\norm{Z}_F} \\
	\eta_{\mcal C1} = \frac{\mid \innerp{X,S}\mid}{1+\norm{X}_F+\norm{S}_F} &\qquad \eta_{\mcal C2} = \frac{\mid \innerp{X,\mcal D^TZ}\mid}{1+\norm{X}_F+\norm{\mcal D^TZ}_F}
\end{array}
\end{align*}
where each column of $X$ is one of the variables $x_{ij}$ (an analogous statement holds for $S$ and $Z$). 

\section{Results}\label{sec:results}

We tested C-SDP\footnote{The code used in this work is available at  \hyperlink{https://github.com/fsbravo/csdp.git}{https://github.com/fsbravo/csdp.git}.} on a variety of sparse graph matching and QAP-type problems. In particular, we used C-SDP to obtain lower and upper bounds on problems from both the QAP and TSP libraries, and to tackle the \emph{assignment problem} in Nuclear Magnetic Resonance Spectroscopy (NMR).

\subsection{QAP and TSP bounds}

SDP relaxations like the one in \cite{Wolkowicz} and C-SDP can yield both a lower bound and an upper bound for QAP-type problems. The latter is given by the objective value of the semidefinite program, while the former is obtained from $\trace{P^TAPB^T}$ after projecting the doubly-stochastic matrix to the set of permutations. In what follows we present both lower and upper bounds for various QAP and TSP problems. 

We compare our results against two other methods: the eigenspace relaxation \cite{deKlerk1}, \cite{deKlerk2}, and the convex-concave approach, PATH \cite{PATH}. To the best of our knowledge, the eigenspace relaxation is the only SDP relaxation for the QAP which can handle larger graphs, although an interior point approach is slow for graphs with more than $50$ nodes. In particular, we compare our lower bounds to the ones produced by this relaxation. Although the eigenspace relaxation also produces a doubly-stochastic matrix, there is no obvious way of recovering the original permutation matrix from this variable (direct projection to the set of permutation matrices, or a Birkhoff-von Neuman decomposition \cite{Dufossé} of the doubly stochastic matrix did not produce meaningful results). We compare our upper bounds against PATH (which produces a permutation matrix). 

For both lower and upper bounds, we compute the gap:
\begin{equation}
\mu = \frac{\mid v^*-v\mid }{v^*} \times 100 \%
\end{equation}
where $v$ is the value of the bounds obtained from the relaxation and $v^*$ is the optimal value for the non-convex problem. The gap can be greater than $100\%$ in the case of the upper bound.

\subsubsection{QAP library problems.}

\begin{figure}[!h]
	\centering
	\includegraphics[width=0.7\columnwidth]{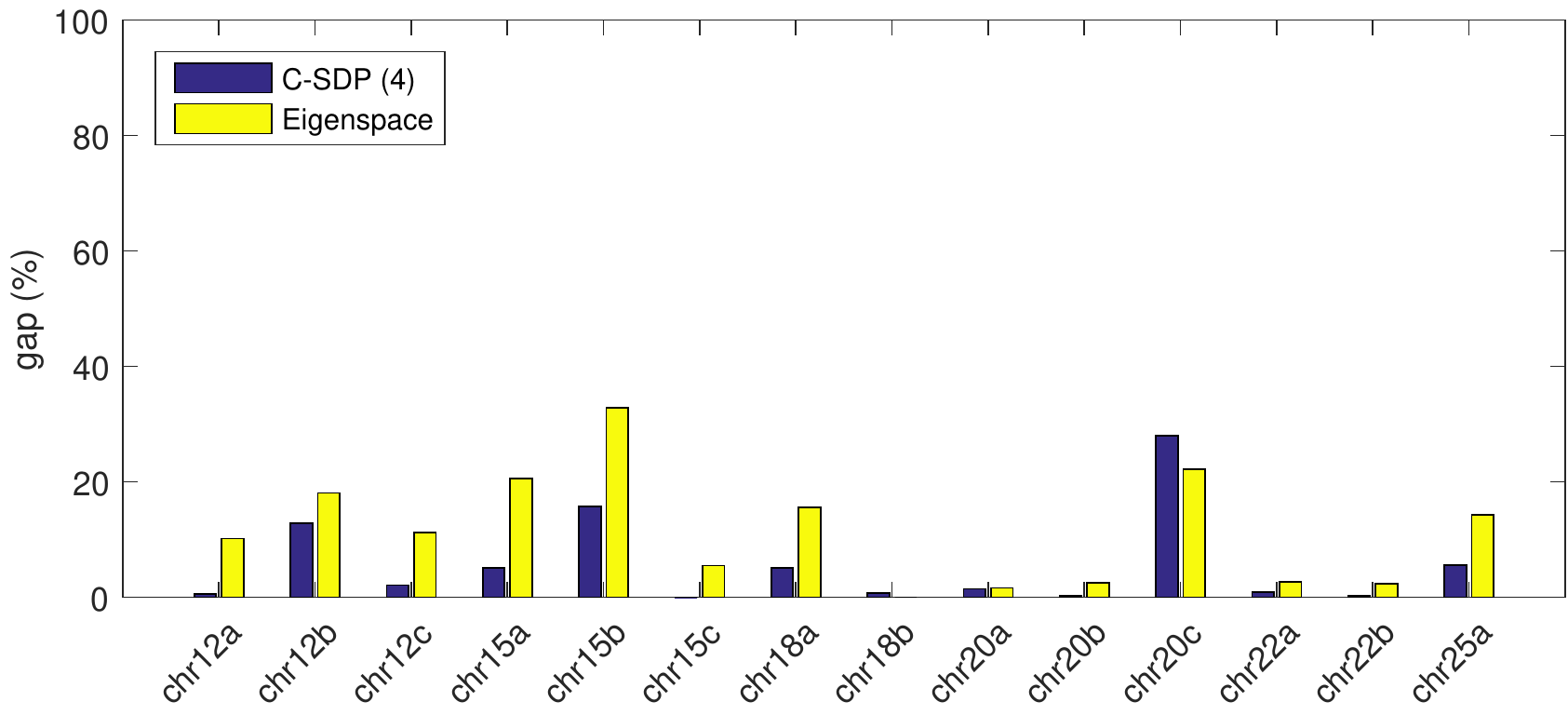}
	\includegraphics[width=0.7\columnwidth]{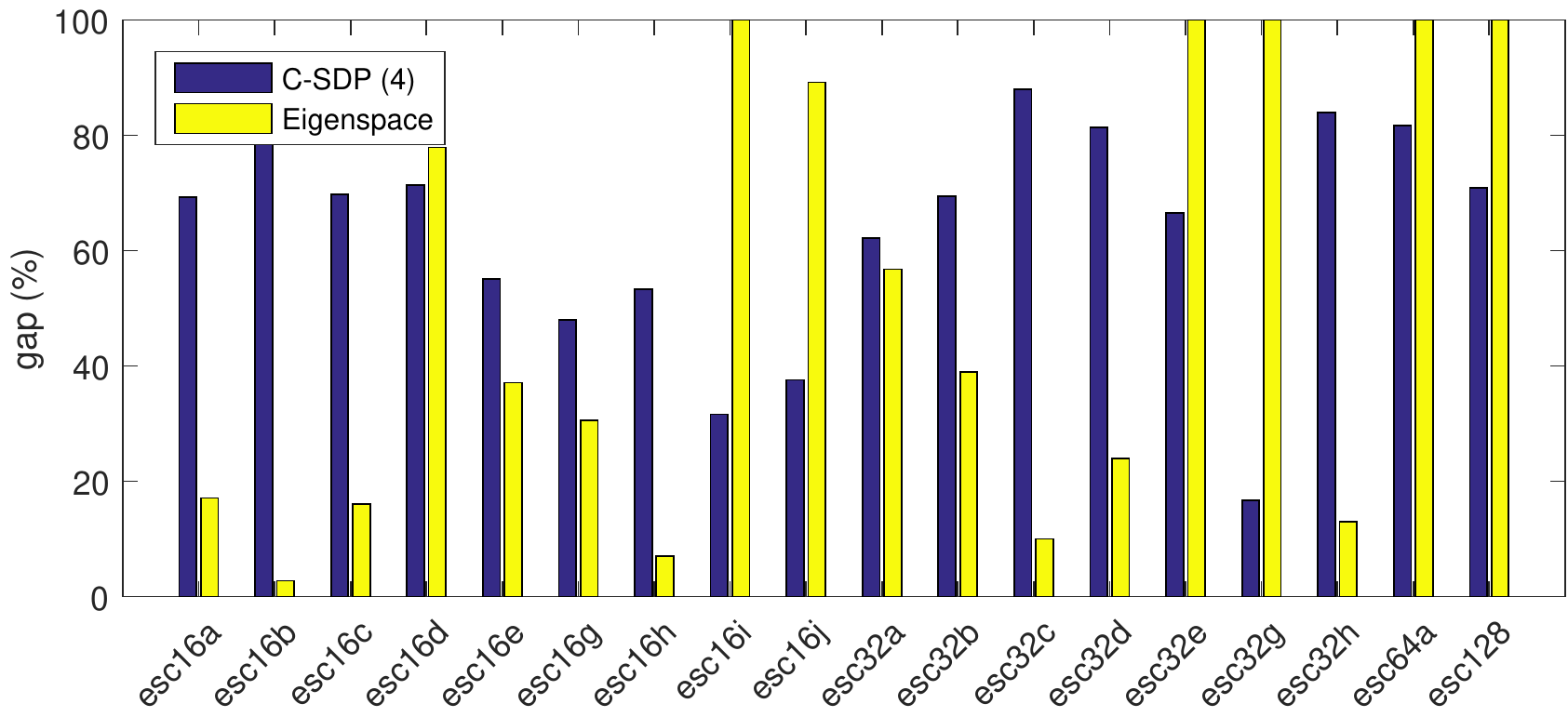}
	\caption{Lower bounds for selected problems in the QAP library. Lower bounds from C-SDP with 4 nodes per variable (blue) and for Eigenspace (yellow) are shown. C-SDP generally shows a smaller gap than Eigenspace for the 'chr' problem family. For the 'esc' problem family, the results are mixed. However, note that out of the 10 problems in which C-SDP underperforms, 7 of these have $B$ with more than 20\% non-zero entries.}\label{fig:qap_lower}
\end{figure}

Figure \ref{fig:qap_lower} shows lower bounds on problems in the QAP library for both C-SDP with cliques of size 4 and Eigenspace (full results are shown in table \ref{tab:qap_lower} in the Appendix). We see that for the 'chr' family of problems \cite{Christofides} C-SDP tends to be significantly better than Eigenspace. However, for the 'esc' problem family \cite{Eschermann90}, the results are divided. Eigenspace performs better than C-SDP in 10 of these problems. In 7 of these, the adjacency matrix, $B$, has more than 20\% non-zero entries. This illustrates a general observed trend, where C-SDP performs best in problems with very sparse $B$.

\begin{figure}[!h]
	\centering
	\includegraphics[width=0.7\columnwidth]{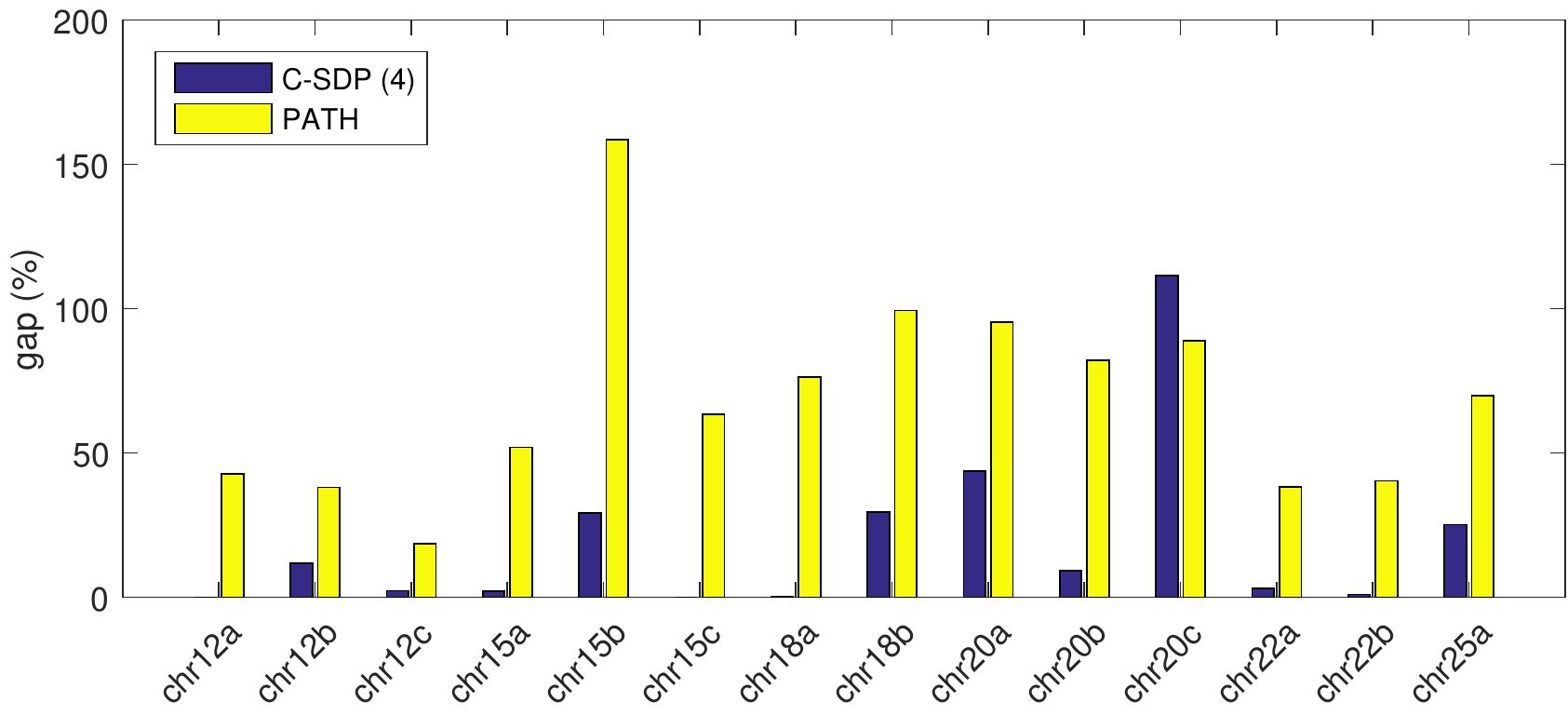}
	\includegraphics[width=0.7\columnwidth]{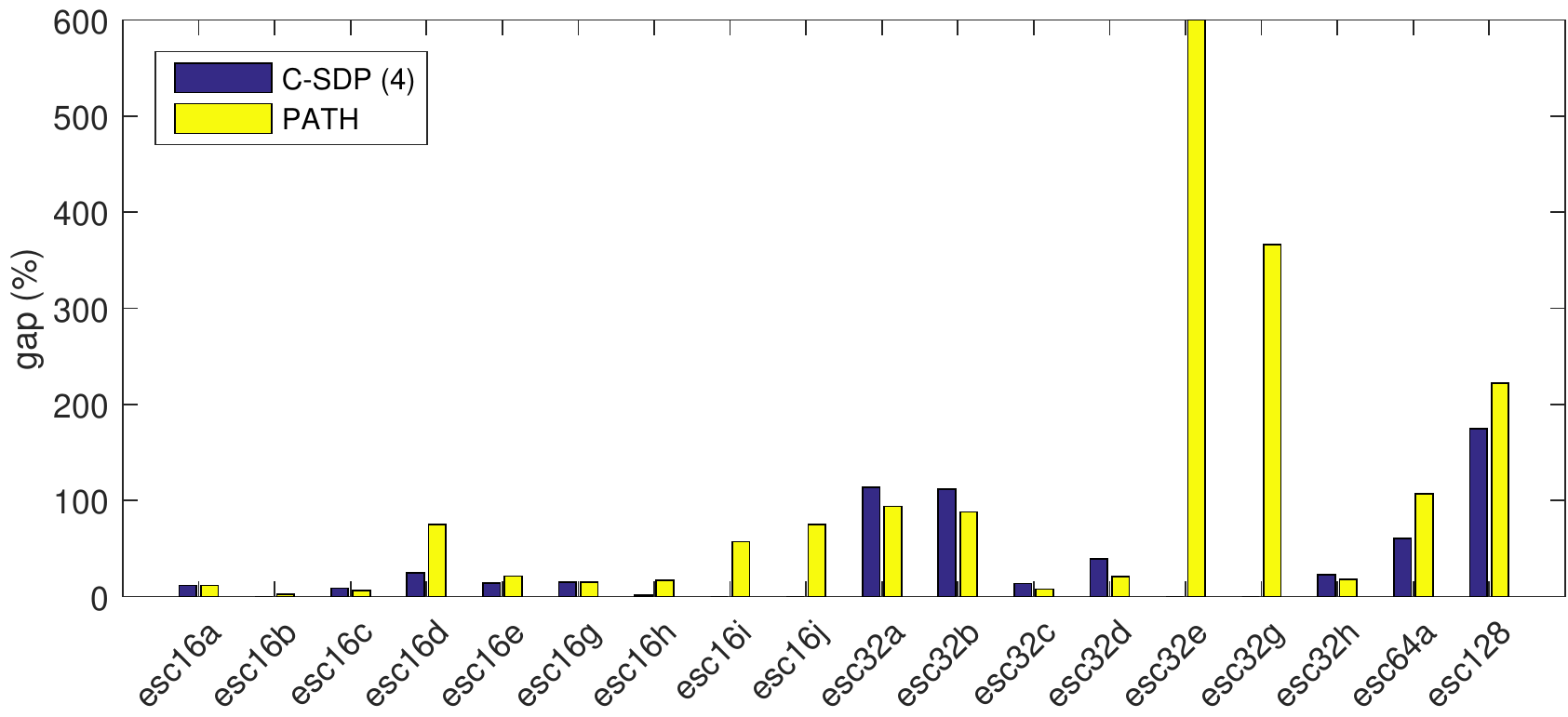}
	\caption{Upper bounds for selected problems in the QAP library. Upper bounds from C-SDP with 4 nodes per variable (blue) and for PATH (yellow) are shown. C-SDP generally shows a smaller gap than PATH.}\label{fig:qap_upper}
\end{figure}

Figure \ref{fig:qap_upper} shows upper bounds on the same problems (full results shown in Table \ref{tab:qap_upper}). This time, a comparison is made with the convex-concave approach, PATH. Generally, we observe that C-SDP produces strong lower bounds on these QAP problems (mostly within 20\% of the optimum). Remarkably, C-SDP achieves the optimum in 7 of the problems. PATH outperforms C-SDP in terms of upper bounds in only 6 of the 32 problems. 

\subsubsection{TSP library problems.}

\begin{figure}[!h]
	\centering
	\includegraphics[width=0.7\columnwidth]{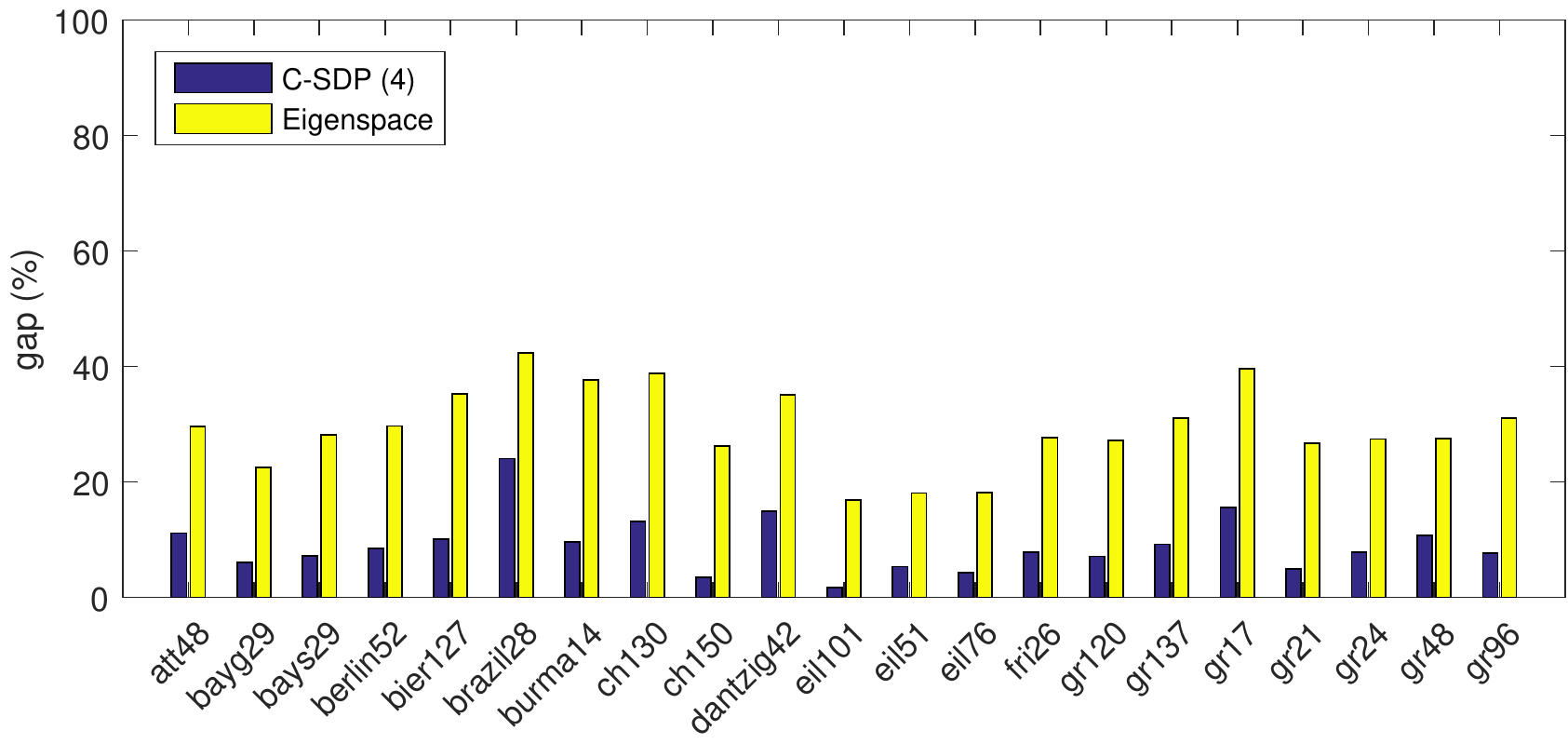}
	\includegraphics[width=0.7\columnwidth]{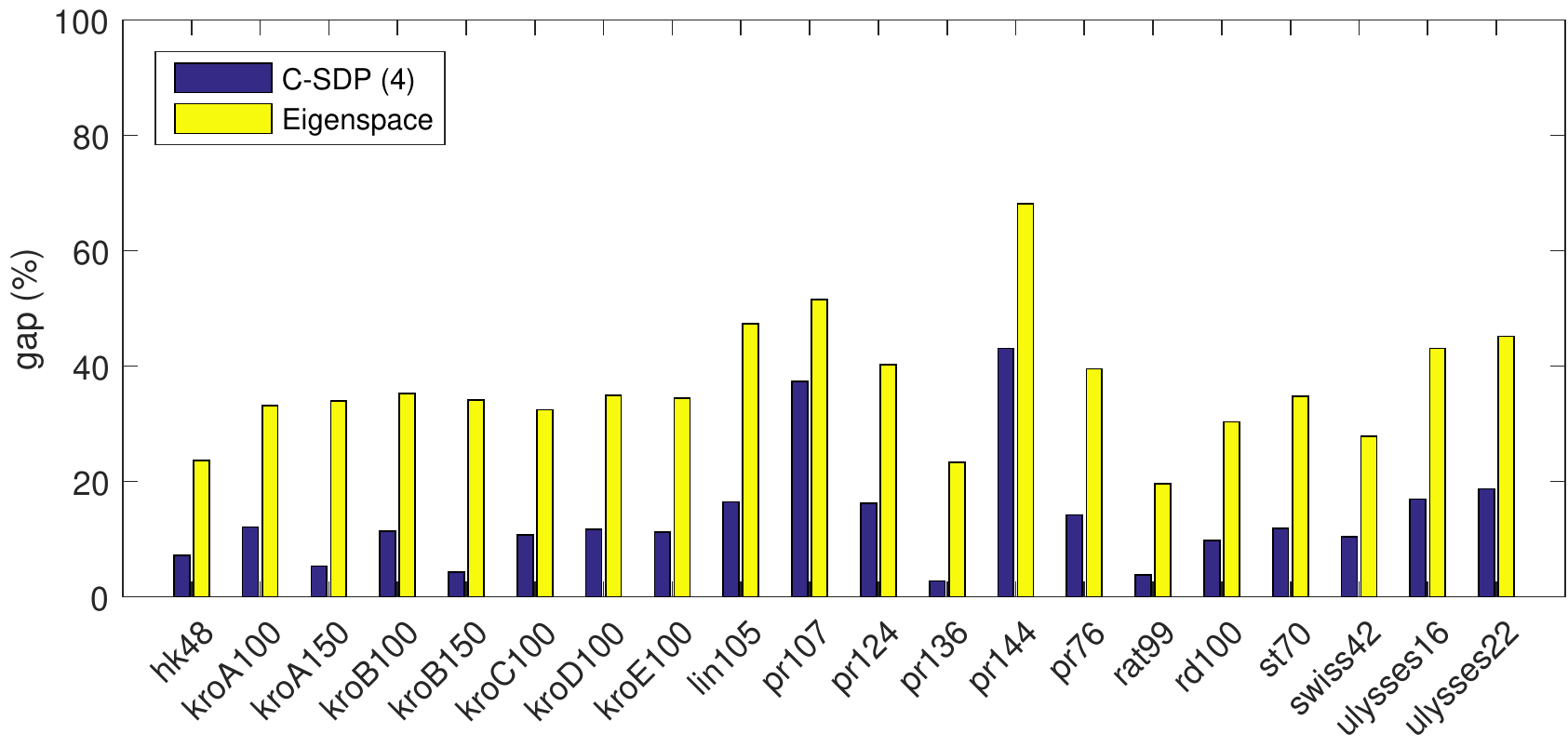}
	\caption{Lower bounds for selected problems from the TSP library ($n\leq 150$). Lower bounds from C-SDP with 4 nodes per variable (blue) and for Eigenspace (yellow) are shown. C-SDP consistently shows smaller gaps than Eigenspace, although Eigenspace can be used to quickly generate a lower bound, since it simplifies to a linear program due to the simple spectrum of B.}\label{fig:tsp_lower}
\end{figure}

\begin{figure}[!h]
	\centering
	\includegraphics[width=0.7\columnwidth]{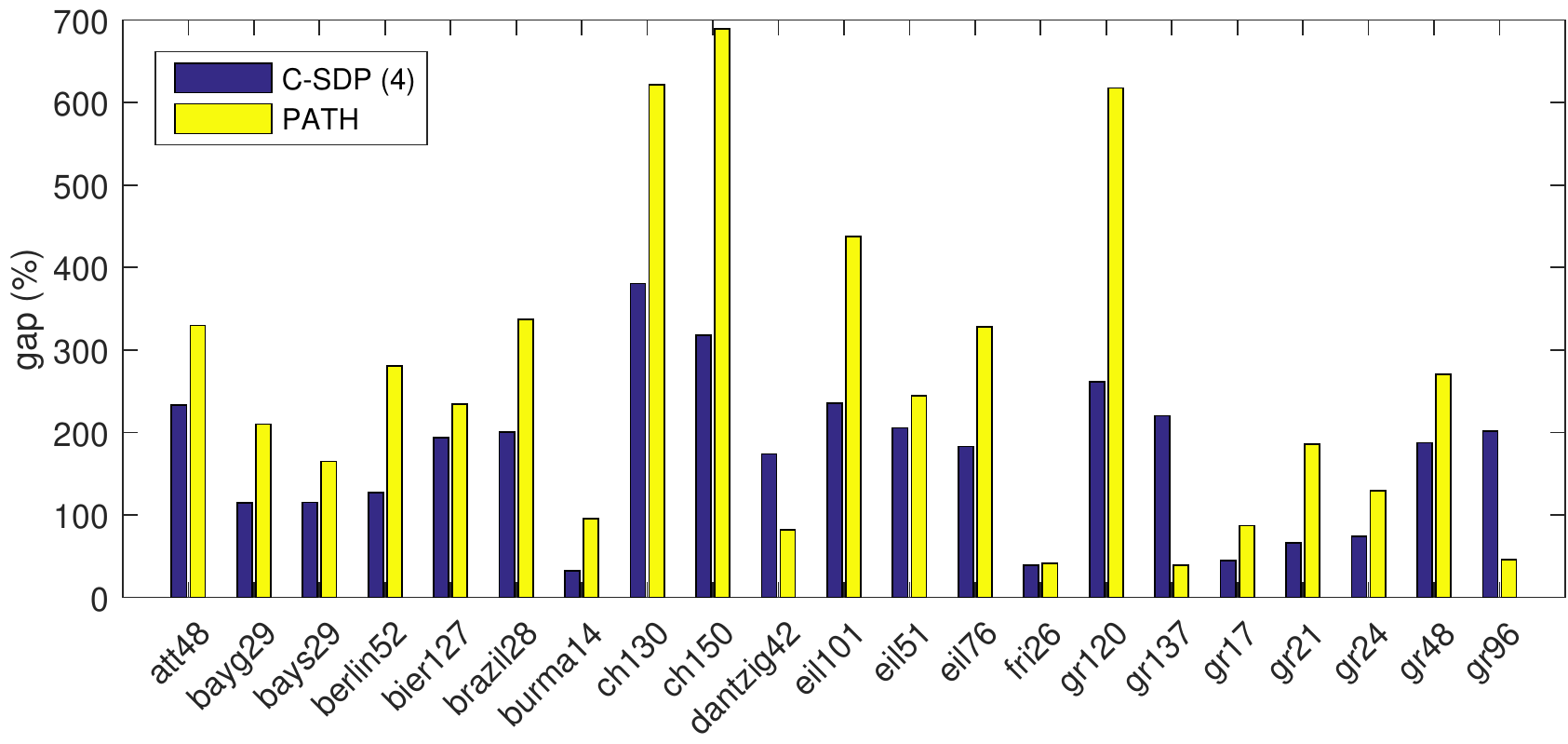}
	\includegraphics[width=0.7\columnwidth]{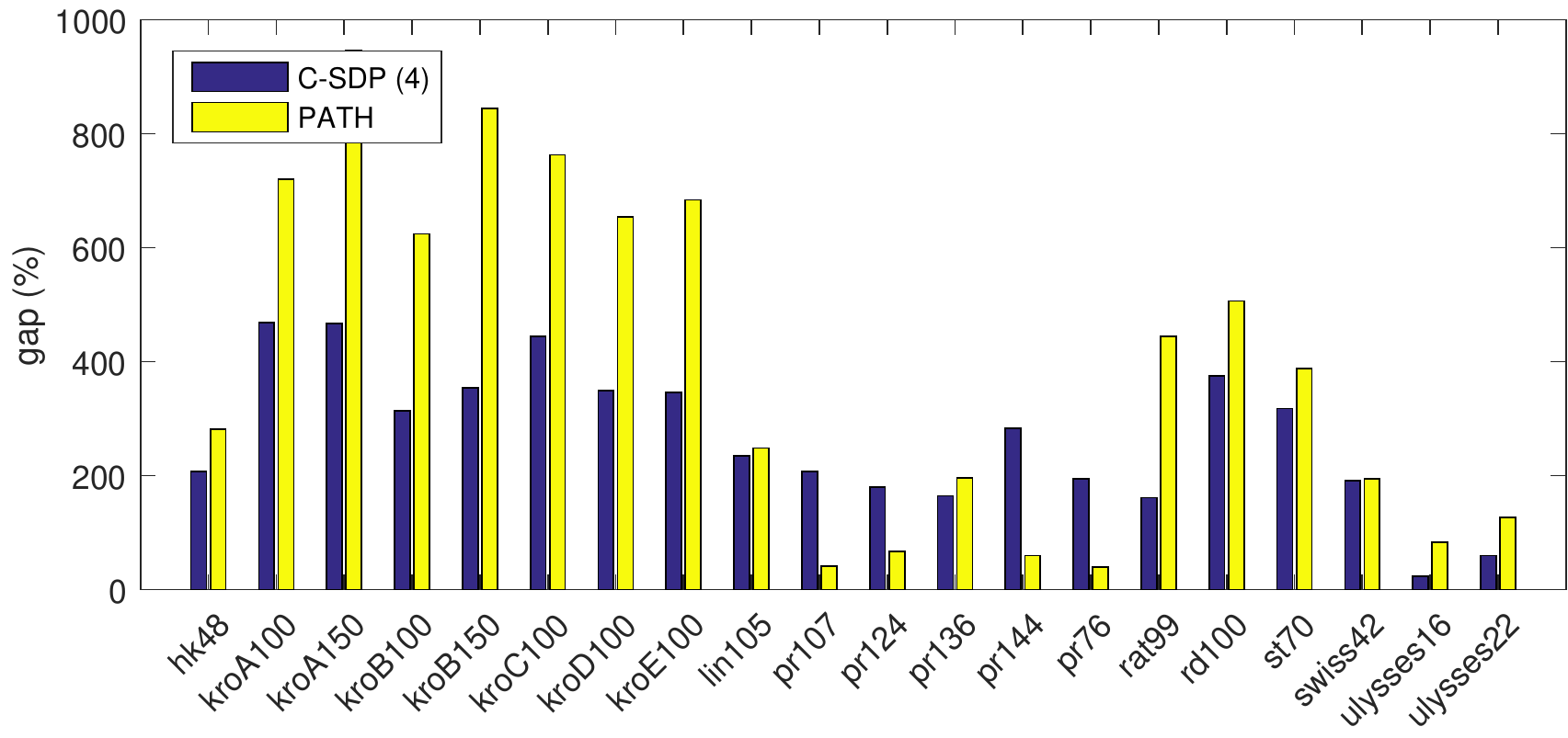}
	\caption{Upper bounds for selected problems in the TSP library ($n\leq 150$). Upper bounds from C-SDP with 4 nodes per variable (blue) and for PATH (yellow) are shown. C-SDP generally shows a smaller gap than PATH.}\label{fig:tsp_upper}
\end{figure}

Figures \ref{fig:tsp_lower} and \ref{fig:tsp_upper} show lower and upper bounds for problems in the TSP library. Again, we verify that C-SDP tends to produce strong lower bounds. Both C-SDP and PATH fail at producing good upper bounds on this class of problems. Tables \ref{tab:tsp_lower} and \ref{tab:tsp_upper} show the detailed results for lower and upper bounds, respectively.

\subsection{Application: NMR assignment}
Nuclear Magnetic Resonance Spectroscopy (NMR) is the go-to tool for structural determination of proteins in solution, \cite{Wuthrich81}. Structural reconstruction in NMR requires accurate geometrical constraints which are derived from the analysis of NMR spectra.

Prior to an NMR experiment, the amino acid sequence (and hence also the atomic composition of the protein) is known. In an experiment, the resonance frequencies of all atoms in the protein are simultaneously measured. In order to use the experimental measurements as constraints on the atoms, measured resonance frequencies have to be assigned to the atoms in the protein, giving rise to the resonance assignment problem. The assignment procedure is typically performed in two steps: 1) Grouping the resonance frequencies from each amino acid into \emph{spin systems}; 2) Assignment of spin systems to the amino acids. One can view spin systems as a vector of resonance frequencies associated with atoms from an amino acid. By assigning each spin system to the correct amino acid, one can then infer the frequencies of the atoms that compose that amino acid.

\begin{figure}
	\centering
	\includegraphics[width=0.7\columnwidth]{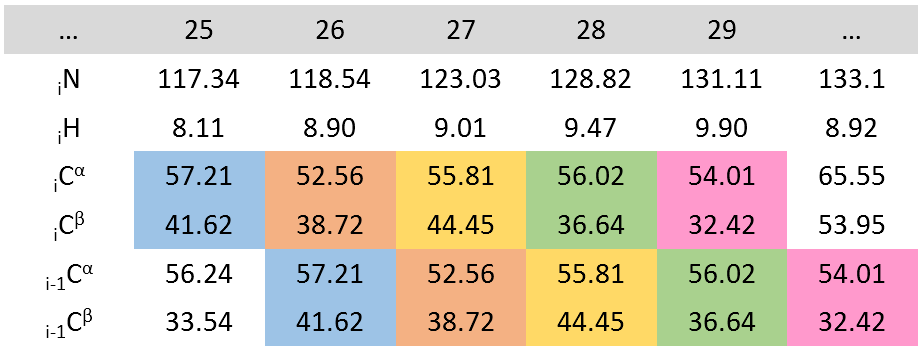}
	\caption{Example sequence of consecutive spin systems built from three NMR spectra for bmr4391 (no noise). Each spin system contains frequency information for $C^\alpha$ and $C^\beta$ in its own amino acid and the preceding amino acid.}\label{fig:ss_example}
\end{figure}

In the following, we formulate the resonance assignment problem as a QAP (the use of the QAP for NMR assignment is not new \cite{Tabu}, \cite{Pistachio}, although to the best of our knowledge, our formulation of the costs and constraints is novel). More precisely, when placed in the correct order, each spin system shares frequency values with its preceding neighbor, up to experimental noise, since some of the atoms of a single amino acid are featured on two adjacent spin systems. This is illustrated in Figure \ref{fig:ss_example}. 

Such a feature can be used to define a distance matrix between each pair of spin systems,
\begin{equation}
A_{ij} = \frac{(s_i(C_i^\alpha)-s_j(C_{j-1}^\alpha))^2}{\sigma_\alpha^2}+\frac{(s_i(C_i^\beta)-s_j(C_{j-1}^\beta))^2}{\sigma_\beta^2}.
\end{equation}
Note that $A_{ij}$ is small if $i$ immediately precedes $j$. One wishes to find the permutation that minimizes this distance along a path of length $n-1$, which should correspond to the best ordering of the spin systems, thus allowing for their assignment to the corresponding amino acids. Alternatively, we choose instead to kernalize this distance matrix by defining
\begin{equation}
\bar A=\exp\left(-\frac{A}{\norm{A}_F}\right)
\end{equation}
where the exponential is applied elementwise \cite{Genton02}. The goal is then to maximize the kernel distance along a path of length $n-1$, amounting to the following problem

\begin{problem}[NMR assignment]
	\begin{align*}
	\max_P & \qquad \trace{\bar APB^TP^T}\\
	\st & \qquad P\in\Perm
	\end{align*}
\end{problem}

where $B$ is the adjacency matrix for the path graph, given by
\begin{equation*}
B = \begin{bmatrix}
0      & 1      & \cdots & \cdots & 0 \\
\vdots & & \ddots &  & \vdots \\
\vdots &        &  &  & 1 \\
0      & \cdots & \cdots & \cdots & 0
\end{bmatrix}
\end{equation*}

This is a QAP, and the sparsity of $B$ allows the use of C-SDP to obtain a doubly stochastic matrix $D$. The matrix $D$ can in turn be projected to the set of permutation matrices to yield a valid assignment of spin systems to amino acids. 

\begin{figure}
	\centering
	\includegraphics[width=0.6\columnwidth]{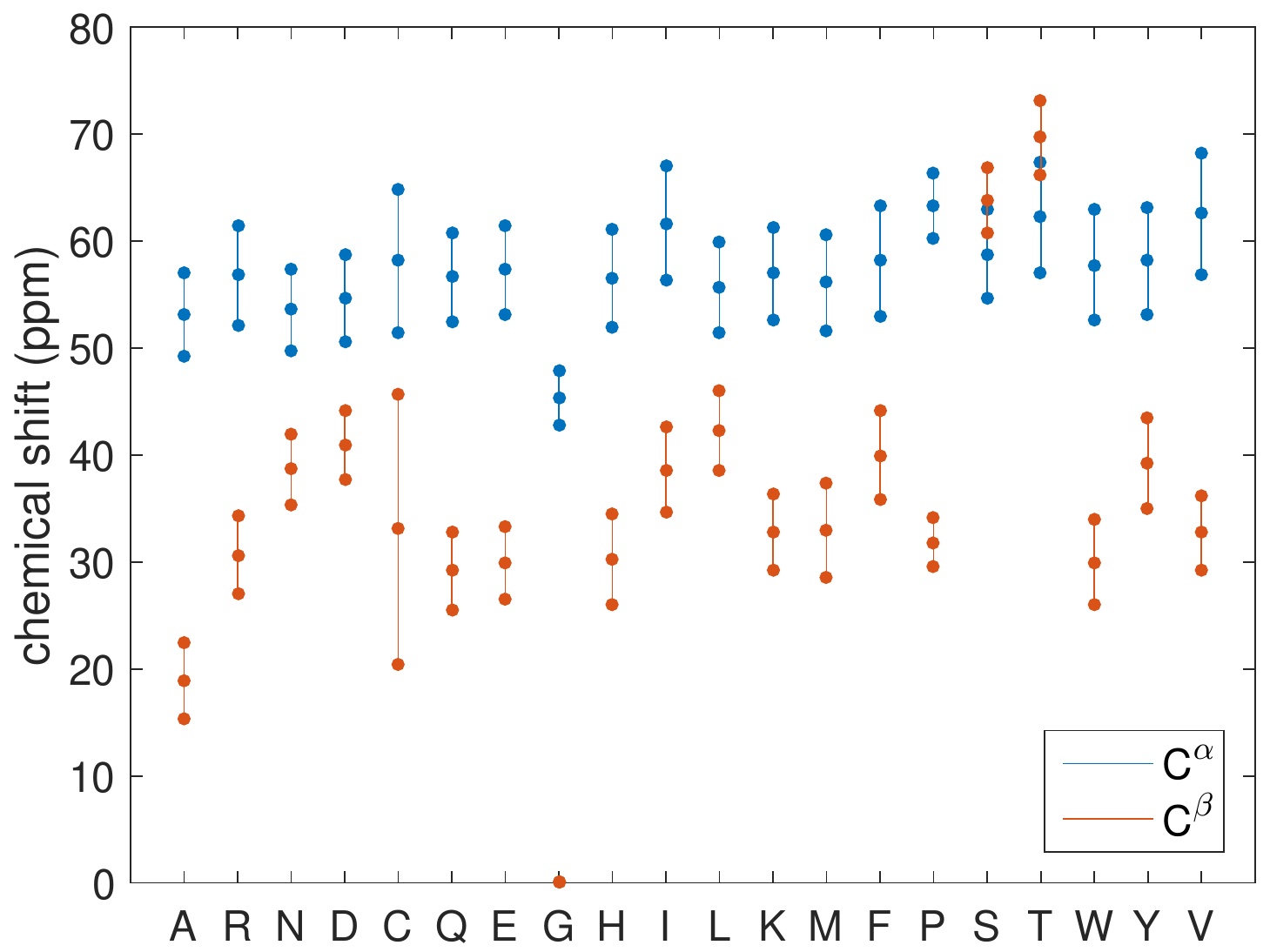}
	\caption{Chemical shift (frequency) distribution of $C^\alpha$ and $C^\beta$ atoms per residue type. The central dots correspond to the mean, and the upper and lower dots are placed two standard deviations away from the mean.}
	\label{fig:cc_distribution}
\end{figure}

Additional information about valid assignments can be included through a term $\trace{W^TP}$ in the QAP cost. In particular, one verifies in practice that the resonance frequencies of the $N$, $H^N$, $C^\alpha$, and $C^\beta$ atoms depend strongly on the type of amino acid, as illustrated in Figure \ref{fig:cc_distribution}. Here we use a hypothesis test perspective to construct $W$. The resonance frequency distributions, conditional on the type of amino acid, were modeled as independent Normal distributions, with mean, $\mu$, and standard deviation, $\sigma$, taken from statistics collected in the Biological Magnetic Resonance Data Bank (BMRB, \cite{BMRB}). Let spin system $s_i$ be defined as the frequencies of the atoms corresponding to that spin system:
\begin{equation}
s_i = \left[\begin{array}{c}
f_{\text{N}_i} \\ f_{\text{H}^N_i} \\ f_{\text{C}^\alpha_i} \\ f_{\text{C}^\beta_i} \\ f_{\text{C}^\alpha_{i-1}} \\ f_{\text{C}^\beta_{i-1}}
\end{array}\right].
\end{equation}
Then for spin system $s_i$ and amino acid $j$ one can compute
\begin{equation*}
z_{ij} = \sum_{k=1}^4 \frac{(s_i(k)-\mu_j(k))^2}{\sigma_k^2} \sim \chi^2_3
\end{equation*}
which follows a chi-square distribution with $3$ degrees of freedom under the distributional assumptions. The test value $z_{ij}$ can be used to determine a p-value under the null hypothesis that spin system $i$ corresponds to residue $j$. Such p-values can be used either as hard constraints on the permutation matrix $P$ (by setting $P(i,j)=0$ if $z_{ij}$ is below a set threshold) or as soft constraints by including an additional term $\trace{W^TP}$ in the cost where $W$ is
\begin{equation*}
W=\exp\left(-\frac{Z}{\norm{Z}_F}\right), 
\end{equation*} with $Z\equiv [z_{ij}]$. The latter approach is adopted in this work. 

Writing the QAP in the form of Problem \ref{prob:1}
\begin{align*}
\min_Q & \qquad \trace{CQ} \\
\st & \qquad Q = \vect{P}\vect{P}^T \\
& \qquad P\in \Perm
\end{align*}
we wish to solve this problem where $C$ is given by
\begin{align*}
C = -B\otimes \exp\left(-\frac{A}{\norm{A}_F}\right)-\frac{1}{\gamma}\text{diag}\left(\exp\left(-\gamma \frac{Z}{\norm{Z}_F}\right)\right).
\end{align*}
Note that the change in the objective consisted only of adding weight terms to the diagonal, using $\gamma$ as a parameter controlling the importance of the statistical information from empirically observed frequencies. A value of $\gamma=0.1$ was used throughout all simulations.

C-SDP was tested on synthetic datasets for benchmarking with cliques of size 2. The dataset was originally described in \cite{Cisa} and consists of a number of proteins for which spin systems are created from the existing assignments in BMRB \cite{BMRB}. We consider only those proteins with 100 amino acids or fewer. Each spin system is constructed by taking the assigned frequencies from the datafile for the base $N$-$H^N$ pair, and the $C^\alpha$ and $C^\beta$ atoms (with the exception of Glycine and Proline). The $C^\alpha$ and $C^\beta$ values from the preceding residue are then added at the end of this vector, and perturbed with additive white gaussian noise with $\sigma=(0.08,0.16)$ (low-noise) or $\sigma=(0.16,0.32)$ (high-noise).

\paragraph{Comparison:} We compare results with other fully automated assignment tools: MARS \cite{Mars}, CISA \cite{Cisa}, and IPASS \cite{Ipass}. In order to compare with a different convex relaxation of the graph matching problem, we consider the doubly-stochastic relaxation (DS), in which we solve the convex problem
\begin{align*}
\min_{D} & \qquad \norm{AD-DB}_F \\
\st & \qquad \trace{K^TD} = 0 \\
& \qquad D \in \text{DS}(n)
\end{align*}
where $K$ is the matrix defined by the entries
\begin{equation*}
K(i,j)=\left\{\begin{array}{cl}
1 & \text{if } z_{ij}<\epsilon \\
0 & \text{otherwise,}
\end{array}\right.
\end{equation*}
for some user-defined threshold $\epsilon$, thus imposing hard constraints on assignments that are statistically unlikely. The value of $\epsilon$ was progressively reduced from $\epsilon=10^{-2}$ until a satisfiable set of constraints was produced.

\paragraph{Evaluation:} Let $N_m$ be the number of spin systems assigned in the BMRB file, $N_{a}$ be the number assigned by the algorithm, and $N_c$ be the number of correctly assigned spin systems. We then define precision$\equiv N_c/N_m$, and recall$\equiv N_c/N_a$. These values are presented in tables \ref{tab:synthetic_low} and \ref{tab:synthetic_high}, below.

\begin{table}[h!]
	\centering
	\caption{ Accuracy of assignment (precision/recall) of various assignment packages as well as the constrained relaxed graph matching (DS) and C-SDP on synthetic spin systems with noise level = (0.08,0.16). Results for MARS \protect\cite{Mars} and CISA taken from \protect\cite{Cisa}. Results for IPASS taken from \protect\cite{Ipass}. }
	\label{tab:synthetic_low}
	\begin{tabular}{lcccccc}
		\hline
		\textbf{Protein ID} & \textbf{Length} & \textbf{MARS} & \textbf{CISA} & \textbf{IPASS} & \textbf{DS} & \textbf{C-SDP} \\ \hline
		bmr4391             & 66              & 100/76        & 97/97         & 93/90          & 85.2/85.2    & 99.1/99.1       \\
		bmr4752             & 68              & 100/97        & 96/94         & 100/94         & 98.5/98.5    & 100/100        \\
		bmr4144             & 78              & 100/91        & 100/99        & 98/85          & 96.4/96.4    & 99.7/99.7       \\
		bmr4579             & 86              & 99/98         & 98/98         & 100/98         & 99.9/99.9    & 100/100        \\
		bmr4316             & 89              & 100/100       & 100/99        & 99/98          & 95.8/95.8    & 98.8/98.8
	\end{tabular}
\end{table}

\begin{table}[h!]
	\centering
	\caption{Accuracy of assignment (precision/recall) of various assignment packages as well as the constrained relaxed graph matching (DS) and C-SDP on synthetic spin systems with noise level = (0.16,0.32). Results for MARS \protect\cite{Mars} and CISA taken from \protect\cite{Cisa}. Results for IPASS taken from \protect\cite{Ipass}. }
	\label{tab:synthetic_high}
	\begin{tabular}{lcccccc}
		\hline
		\textbf{Protein ID} & \textbf{Length} & \textbf{MARS} & \textbf{CISA} & \textbf{IPASS} & \textbf{DS} & \textbf{C-SDP} \\ \hline
		bmr4391             & 66              & 100/75        & 91/91         & 93/90          & 85.5/85.5    & 100/100        \\
		bmr4752             & 68              & 100/97        & 90/88         & 100/94         & 87.8/87.8    & 99.4/99.4       \\
		bmr4144             & 78              & 100/69        & 100/99        & 98/85          & 85.6/85.6    & 96.4/96.4       \\
		bmr4579             & 86              & 96/90         & 80/80         & 100/98         & 89.6/89.6    & 99.6/99.6       \\
		bmr4316             & 89              & 99/91       & 83/83        & 99/98          & 95.1/95.1    & 97.8/97.8
	\end{tabular}
\end{table}

C-SDP outperforms other methods in terms of both precision and recall, on average. The doubly stochastic relaxation also compares well in the low-noise scenario, but performs poorly in the high-noise setting. 

Note that since C-SDP always produces a full permutation matrix, it assigns all spin systems, such that $N_a=N_m$, where $N_m$ is assumed to be the number of assignable spin systems in the protein. As a result, precision and recall values for C-SDP are equal. The number of spin systems may be smaller than the number of amino acids, as some amino acids such as Proline do not contribute spin systems, in which case token spin systems are used in their stead.
\section{Conclusion}

This work presented a new semidefinite programming (SDP) relaxation, C-SDP, for the quadratic assignment problem (QAP) in which one of the matrices is sparse. A convergent ADMM formulation was developed, which exploits the natural three-block structure of the dual problem, allowing a highly parallelizable solution where the most expensive step per iteration is a projection of a matrix of size $\mcal{O}(n)$ to the positive semidefinite cone.

The performance of C-SDP was evaluated on problems from the QAP and TSP libraries, where we found it produces better lower bounds than comparable SDP relaxations \cite{deKlerk2} as well as competitive upper bounds (after projecting the solution to the set of permutation matrices), compared to a popular local method \cite{PATH}.

An application to the NMR assignment problem was also described, which can be formulated as a sparse QAP. Preliminary results on proteins from a standard synthetic dataset showed that C-SDP results in a better assignment compared to popular fully automated assignment tools in recent literature.

\section*{\large Acknowledgements}\small
	The authors would like to thank David Cowburn for useful discussions on NMR spectroscopy, and Amir Ali Ahmadi, for suggestions on tightening the SDP relaxations presented here.
	
	The authors were partially supported by Award Number R01GM090200
	from the NIGMS, FA9550-12-1-0317 from AFOSR, the Simons Investigator
	Award and the Simons Collaboration on Algorithms and Geometry from Simons
	Foundation, and the Moore Foundation Data-Driven Discovery Investigator
	Award.

\bibliographystyle{spmpsci} 
\bibliography{references} 


\newpage
\section*{Appendix - Full Tabulated Results}
\begin{table}[!h]\small
	\centering
	\caption{Comparison between lower bounds given by the C-SDP and Eigenspace relaxations on selected problems from the QAP library with (relatively) sparse B. C-SDP instances were ran for 1000 ADMM iterations on 20 processors.}
	\label{tab:qap_lower}
	\begin{tabular}{lccccc}
		\hline
		\textbf{Problem} & \textbf{Optimal} & \textbf{\begin{tabular}[c]{@{}c@{}}C-SDP\\ (k=2)\\ Gap (\%)\end{tabular}} & \textbf{\begin{tabular}[c]{@{}c@{}}C-SDP\\ (k=3)\\ Gap (\%)\end{tabular}} & \textbf{\begin{tabular}[c]{@{}c@{}}C-SDP\\ (k=4)\\ Gap (\%)\end{tabular}} & \textbf{\begin{tabular}[c]{@{}c@{}}Eigen-\\ space\\ Gap (\%)\end{tabular}} \\ \hline
		chr12a  & 9552    & 9.7          & 2.7          & \textbf{0.6}  & 10.2          \\
		chr12b  & 9742    & 26.3         & 19.3         & \textbf{12.8} & 18.1          \\
		chr12c  & 11156   & 10.2         & 3.2          & \textbf{2.1}  & 11.3          \\
		chr15a  & 9896    & 13.1         & 7.0          & \textbf{5.2}  & 20.6          \\
		chr15b  & 7990    & 35.7         & 25.2         & \textbf{15.7} & 32.8          \\
		chr15c  & 9504    & 0.1          & \textbf{0.0} & \textbf{0.0}  & 5.5           \\
		chr18a  & 11098   & 12.6         & \textbf{4.6} & 5.2           & 15.5          \\
		chr18b  & 1534    & \textbf{0.0} & \textbf{0.0} & 0.7           & \textbf{0.0}  \\
		chr20a  & 2192    & 1.6          & 1.6          & \textbf{1.5}  & 1.6           \\
		chr20b  & 2298    & 2.8          & 0.3          & \textbf{0.3}  & 2.5           \\
		chr20c  & 14142   & 37.0         & 31.9         & 28.0          & \textbf{22.3} \\
		chr22a  & 6156    & 2.6          & 1.3          & \textbf{0.9}  & 2.7           \\
		chr22b  & 6194    & 1.4          & 0.5          & \textbf{0.2}  & 2.4           \\
		chr25a  & 3796    & 13.8         & 7.5          & \textbf{5.6}  & 14.3          \\
		esc16a  & 68      & 100.0        & 80.4         & 69.2          & \textbf{17.1} \\
		esc16b  & 292     & 100.0        & 88.4         & 85.7          & \textbf{2.7}  \\
		esc16c  & 160     & 100.0        & 79.4         & 69.8          & \textbf{16.1} \\
		esc16d  & 16      & 100.0        & 72.1         & \textbf{71.4} & 77.9          \\
		esc16e  & 28      & 100.0        & 78.2         & 55.0          & \textbf{37.1} \\
		esc16g  & 26      & 100.0        & 67.3         & 48.0          & \textbf{30.6} \\
		esc16h  & 996     & 66.3         & 57.0         & 53.3          & \textbf{7.0}  \\
		esc16i  & 14      & 100.0        & \textbf{19.1}& 31.6          & 100.0         \\
		esc16j  & 8       & 100.0        & 53.0         & \textbf{37.6} & 89.2          \\
		esc32a  & 130     & 100.0        & 74.6         & 62.1          & \textbf{56.8} \\
		esc32b  & 168     & 100.0        & 76.1         & 69.5          & \textbf{39.0} \\
		esc32c  & 642     & 100.0        & 91.0         & 88.0          & \textbf{10.0} \\
		esc32d  & 200     & 100.0        & 88.0         & 81.3          & \textbf{24.0} \\
		esc32e  & 2       & 100.0        & 66.6         & \textbf{66.5} & 100.0         \\
		esc32g  & 6       & 100.0        & 41.0         & \textbf{16.7} & 100.0         \\
		esc32h  & 438     & 100.0        & 88.2         & 84.0          & \textbf{13.0} \\
		esc64a  & 116     & 99.9         & 86.0         & 81.7          & -             \\
		esc128  & 64      & 99.6         & 75.8         & 71.0          & -             \\
		ste36a  & 96772   & 46.9         & 42.1         & 40.2          & -             \\
		ste36b  & 58537   & 70.1         & 67.0         & 65.1          & -             \\
		ste36c  & 108159  & 37.7         & 35.2         & 34.2          & -
	\end{tabular}
\end{table}

\begin{table}[!h]\small
	\centering
	\caption{Comparison between upper bounds given by C-SDP and PATH relaxations on elected problems from the QAP library with (relatively) sparse B. C-SDP instances were ran for 1000 ADMM iterations on 20 processors.}
	\label{tab:qap_upper}
	\begin{tabular}{lccccc}
		\hline
		\textbf{Problem} & \textbf{Optimal} & \textbf{\begin{tabular}[c]{@{}c@{}}C-SDP\\ (k=2)\\ Gap (\%)\end{tabular}} & \textbf{\begin{tabular}[c]{@{}c@{}}C-SDP\\ (k=3)\\ Gap (\%)\end{tabular}} & \textbf{\begin{tabular}[c]{@{}c@{}}C-SDP\\ (k=4)\\ Gap (\%)\end{tabular}} & \textbf{\begin{tabular}[c]{@{}c@{}}PATH\\ \\ Gap (\%)\end{tabular}} \\ \hline
		chr12a  & 9552    & 34.5  & 6.0           & \textbf{0.0}   & 42.7           \\
		chr12b  & 9742    & 38.9  & 25.4          & \textbf{11.9}  & 38.1           \\
		chr12c  & 11156   & 5.8   & \textbf{2.3}  & \textbf{2.3}   & 18.6           \\
		chr15a  & 9896    & \textbf{2.1}   & \textbf{2.1}  & \textbf{2.1}   & 52.0  \\
		chr15b  & 7990    & \textbf{26.3}  & 34.5          & 29.2           & 158.6 \\
		chr15c  & 9504    & \textbf{0.0}   & \textbf{0.0}  & \textbf{0.0}   & 63.3  \\
		chr18a  & 11098   & 69.8  & \textbf{0.2}  & \textbf{0.2}   & 76.3           \\
		chr18b  & 1534    & \textbf{8.9}   & 22.9 & 29.5           & 99.3           \\
		chr20a  & 2192    & 122.5 & 76.1          & \textbf{43.8}  & 95.4           \\
		chr20b  & 2298    & 62.9  & \textbf{9.3}  & \textbf{9.3}   & 82.2           \\
		chr20c  & 14142   & 173.0 & 100.1         & 111.5          & \textbf{88.9}  \\
		chr22a  & 6156    & 17.2  & 7.6           & \textbf{3.0}   & 38.3           \\
		chr22b  & 6194    & 7.3   & 2.3           & \textbf{1.0}   & 40.4           \\
		chr25a  & 3796    & 107.0 & 49.2          & \textbf{25.2}  & 69.9           \\
		esc16a  & 68      & \textbf{8.8}   & 11.8 & 11.8           & 11.8           \\
		esc16b  & 292     & \textbf{0.0}   & 0.7  & \textbf{0.0}   & 2.7            \\
		esc16c  & 160     & \textbf{5.0}   & 7.5  & 8.7            & 6.3            \\
		esc16d  & 16      & \textbf{12.5} & 50.0  & 25.0           & 75.0           \\
		esc16e  & 28      & 14.3  & \textbf{7.1}  & 14.3           & 21.4           \\
		esc16g  & 26      & 7.7   & \textbf{0.0}  & 15.4           & 15.4           \\
		esc16h  & 996     & 1.6   & \textbf{0.0}  & 1.6            & 16.9           \\
		esc16i  & 14      & \textbf{0.0}   & \textbf{0.0}  & \textbf{0.0}   & 57.1  \\
		esc16j  & 8       & \textbf{0.0}   & \textbf{0.0}  & \textbf{0.0}   & 75.0  \\
		esc32a  & 130     & 115.4 & 124.6         & 113.8          & \textbf{93.8}  \\
		esc32b  & 168     & 109.5 & 114.3         & 111.9          & \textbf{88.1}  \\
		esc32c  & 642     & 12.8  & 15.9          & 13.7           & \textbf{7.8}   \\
		esc32d  & 200     & 38.0  & 36.0          & 39.0           & \textbf{21.0}  \\
		esc32e  & 2       & \textbf{0.0}   & \textbf{0.0}  & \textbf{0.0}   & 600.0 \\
		esc32g  & 6       & \textbf{0.0}   & \textbf{0.0}  & \textbf{0.0}   & 366.7 \\
		esc32h  & 438     & 24.7  & 26.9          & 22.8           & \textbf{18.3}  \\
		esc64a  & 116     & 60.3  & \textbf{53.4} & 60.3           & 106.9          \\
		esc128  & 64      & 250.0 & 206.3         & \textbf{175.0} & 221.9          \\
		ste36a  & 96772   & \textbf{70.2}  & 74.7 & 74.2           & 76.3           \\
		ste36b  & 58537   & 188.8 & 204.3         & 211.9          & \textbf{158.6} \\
		ste36c  & 108159  & 66.0  & \textbf{62.8} & 63.7           & 83.2
	\end{tabular}
\end{table}

\begin{table}[!h]\small
	\centering
	\caption{Comparison between lower bounds given by the C-SDP and Eigenspace relaxations on problems from the TSP library (with $n\leq 150$). C-SDP instances were ran for 1000 ADMM iterations on 20 processors.}
	\label{tab:tsp_lower}
	\begin{tabular}{lccccc}
		\hline
		\textbf{Problem} & \textbf{Optimal} & \textbf{\begin{tabular}[c]{@{}l@{}}C-SDP \\ (k=2)\\ Gap (\%)\end{tabular}} & \textbf{\begin{tabular}[c]{@{}l@{}}C-SDP \\ (k=3)\\ Gap (\%)\end{tabular}} & \textbf{\begin{tabular}[c]{@{}l@{}}C-SDP \\ (k=4)\\ Gap (\%)\end{tabular}} & \textbf{\begin{tabular}[c]{@{}l@{}}Eigen-\\ space\\ Gap (\%)\end{tabular}} \\ \hline
		att48     & 10628   & 20.6  & \textbf{10.9} & 11.2          & 29.6  \\
		bayg29    & 1610    & 10.7  & 6.5           & \textbf{6.0}  & 22.5  \\
		bays29    & 2020    & 12.8  & \textbf{6.9}  & 7.2           & 28.1  \\
		berlin52  & 7542    & 16.5  & 9.9           & \textbf{8.4}  & 29.7  \\
		bier127   & 118282  & 19.5  & 10.4          & \textbf{10.1} & 35.2  \\
		brazil58  & 25395   & 34.6  & 24.1          & \textbf{24.0} & 42.4  \\
		burma14   & 3323    & 16.1  & 10.7          & \textbf{9.7}  & 37.7  \\
		ch130     & 6110    & 26.9  & \textbf{12.9} & 13.1          & 38.8  \\
		ch150     & 6528    & 11.3  & \textbf{1.1}  & 3.5           & 26.3  \\
		dantzig42 & 699     & 23.6  & \textbf{14.9} & 15.0          & 35.1  \\
		eil101    & 629     & 6.9   & \textbf{1.6}  & 1.7           & 16.9  \\
		eil51     & 426     & 10.9  & \textbf{4.6}  & 5.4           & 18.1  \\
		eil76     & 538     & 8.8   & \textbf{3.6}  & 4.3           & 18.2  \\
		fri26     & 937     & 12.2  & 9.4           & \textbf{7.8}  & 27.6  \\
		gr120     & 6942    & 15.1  & \textbf{7.1}  & \textbf{7.1}  & 27.2  \\
		gr137     & 69853   & 14.9  & \textbf{8.0}  & 9.2           & 31.0  \\
		gr17      & 2085    & 21.6  & 20.2          & \textbf{15.6} & 39.7  \\
		gr21      & 2707    & 10.8  & \textbf{4.0}  & 5.0           & 26.7  \\
		gr24      & 1272    & 17.3  & 8.7           & \textbf{7.9}  & 27.4  \\
		gr48      & 5046    & 18.0  & \textbf{10.4} & 10.8          & 27.5  \\
		gr96      & 55209   & 15.0  & 8.4           & \textbf{7.7}  & 31.0  \\
		hk48      & 11461   & 13.9  & 7.5           & \textbf{7.2}  & 23.6  \\
		kroA100   & 21282   & 19.4  & 12.7          & \textbf{12.1} & 33.2  \\
		kroA150   & 26524   & 16.6  & \textbf{5.2}  & 5.3           & 34.0  \\
		kroB100   & 22141   & 23.7  & 12.9          & \textbf{11.4} & 35.3  \\
		kroB150   & 26130   & 18.6  & 4.4           & \textbf{4.3}  & 34.1  \\
		kroC100   & 20749   & 19.3  & 10.9          & \textbf{10.8} & 32.4  \\
		kroD100   & 21294   & 22.2  & 12.1          & \textbf{11.7} & 34.9  \\
		kroE100   & 22068   & 23.9  & 12.5          & \textbf{11.2} & 34.4  \\
		lin105    & 14379   & 35.2  & 17.3          & \textbf{16.5} & 47.4  \\
		pr107     & 44303   & 40.2  & 38.9          & \textbf{37.4} & 51.5  \\
		pr124     & 59030   & 25.3  & \textbf{14.0} & 16.2          & 40.3  \\
		pr136     & 96772   & \textbf{2.1}   & 2.6  & 2.7           & 23.3  \\
		pr144     & 58537   & 54.7  & \textbf{40.7} & 43.1          & 68.2  \\
		pr76      & 108159  & 28.3  & 15.3          & \textbf{14.2} & 39.5  \\
		rat99     & 1211    & 10.4  & \textbf{2.8}  & 3.8           & 19.7  \\
		rd100     & 7910    & 16.9  & 10.7          & \textbf{9.8}  & 30.3  \\
		st70      & 675     & 22.1  & 13.6          & \textbf{11.9} & 34.8  \\
		swiss42   & 1273    & 20.6  & 10.7          & \textbf{10.4} & 27.9  \\
		ulysses16 & 6859    & 25.2  & 17.1          & \textbf{17.0} & 43.1  \\
		ulysses22 & 7013    & 28.5  & 19.7          & \textbf{18.7} & 45.2
	\end{tabular}
\end{table}

\begin{table}[!h]\small
	\centering
	\caption{Comparison between upper bounds given by C-SDP and PATH relaxations on problems from the TSP library (with $n\leq 150$). C-SDP instances were ran for 1000 ADMM iterations on 20 processors.}
	\label{tab:tsp_upper}
	\begin{tabular}{lccccc}
		\hline
		\textbf{Problem} & \textbf{Optimal} & \multicolumn{1}{c}{\textbf{\begin{tabular}[c]{@{}c@{}}C-SDP\\ (k=2)\\ Gap (\%)\end{tabular}}} & \multicolumn{1}{c}{\textbf{\begin{tabular}[c]{@{}c@{}}C-SDP\\ (k=3)\\ Gap (\%)\end{tabular}}} & \multicolumn{1}{c}{\textbf{\begin{tabular}[c]{@{}c@{}}C-SDP\\ (k=4)\\ Gap (\%)\end{tabular}}} & \multicolumn{1}{c}{\textbf{\begin{tabular}[c]{@{}c@{}}PATH\\ \\ Gap (\%)\end{tabular}}} \\ \hline
		att48     & 10628   & \textbf{213.0} & 236.5    & 233.6          & 329.8         \\
		bayg29    & 1610    & \textbf{114.3} & 115.8    & \textbf{114.3} & 210.1         \\
		bays29    & 2020    & \textbf{107.6} & 118.3    & 115.4          & 164.8         \\
		berlin52  & 7542    & 175.0    & \textbf{127.2} & \textbf{127.2} & 280.6         \\
		bier127   & 118282  & 216.4    & \textbf{193.8} & \textbf{193.8} & 234.2         \\
		brazil58  & 25395   & 248.0    & \textbf{200.8} & \textbf{200.8} & 337.0         \\
		burma14   & 3323    & \textbf{24.6}  & 28.4     & 32.3           & 95.5          \\
		ch130     & 6110    & \textbf{352.4} & 380.6    & 380.6          & 621.3         \\
		ch150     & 6528    & 346.9    & \textbf{318.2} & \textbf{318.2} & 689.3         \\
		dantzig42 & 699     & 193.1    & 174.0          & 174.0          & \textbf{82.0} \\
		eil101    & 629     & \textbf{227.3} & 235.3    & 235.3          & 437.7         \\
		eil51     & 426     & \textbf{203.6} & 205.4    & 205.5          & 244.4         \\
		eil76     & 538     & 282.9    & \textbf{183.0} & \textbf{183.0} & 328.2         \\
		fri26     & 937     & 91.6     & \textbf{39.4}  & \textbf{39.4}  & 41.6          \\
		gr120     & 6942    & 445.2    & \textbf{261.6} & \textbf{261.6} & 617.6         \\
		gr137     & 69853   & 264.6    & 220.3          & 220.3          & \textbf{38.9} \\
		gr17      & 2085    & 46.8     & \textbf{32.4}  & 44.9           & 86.9          \\
		gr21      & 2707    & 94.5     & 69.7           & \textbf{66.3}  & 185.7         \\
		gr24      & 1272    & 89.2     & 86.2           & \textbf{73.9}  & 129.4         \\
		gr48      & 5046    & 210.2    & \textbf{187.4} & \textbf{187.4} & 270.4         \\
		gr96      & 55209   & 228.9    & 201.7          & 201.7          & \textbf{46.0} \\
		hk48      & 11461   & 222.4    & \textbf{207.7} & \textbf{207.7} & 281.6         \\
		kroA100   & 21282   & 469.6    & \textbf{469.0} & \textbf{469.0} & 720.2         \\
		kroA150   & 26524   & \textbf{411.0} & 467.4    & 467.4          & 945.8         \\
		kroB100   & 22141   & 411.9    & \textbf{313.6} & \textbf{313.6} & 624.2         \\
		kroB150   & 26130   & 417.3    & \textbf{353.7} & \textbf{353.7} & 844.7         \\
		kroC100   & 20749   & 507.4    & \textbf{445.1} & \textbf{445.1} & 763.0         \\
		kroD100   & 21294   & 504.2    & \textbf{349.8} & \textbf{349.8} & 654.4         \\
		kroE100   & 22068   & 489.5    & \textbf{346.3} & \textbf{346.3} & 684.2         \\
		lin105    & 14379   & 303.1    & \textbf{234.8} & \textbf{234.8} & 248.4         \\
		pr107     & 44303   & 181.5    & 207.9          & 207.9          & \textbf{41.6} \\
		pr124     & 59030   & 293.8    & 180.2          & 180.2          & \textbf{67.6} \\
		pr136     & 96772   & 325.5    & \textbf{164.7} & \textbf{164.7} & 196.6         \\
		pr144     & 58537   & 255.0    & 283.7          & 283.7          & \textbf{59.8} \\
		pr76      & 108159  & 192.2    & 194.0          & 194.0          & \textbf{39.4} \\
		rat99     & 1211    & 236.4    & \textbf{161.5} & \textbf{161.5} & 444.1         \\
		rd100     & 7910    & 438.4    & \textbf{375.3} & \textbf{375.3} & 506.5         \\
		st70      & 675     & \textbf{300.9} & 320.0    & 317.9          & 387.9         \\
		swiss42   & 1273    & \textbf{163.2} & 190.4    & 190.8          & 194.0         \\
		ulysses16 & 6859    & 23.6     & \textbf{20.2}  & 23.2           & 82.7          \\
		ulysses22 & 7013    & 64.5     & \textbf{57.0}  & 59.7           & 126.3
	\end{tabular}
\end{table}

\newpage

\begin{figure}
	\centering
	\includegraphics[width=0.75\columnwidth]{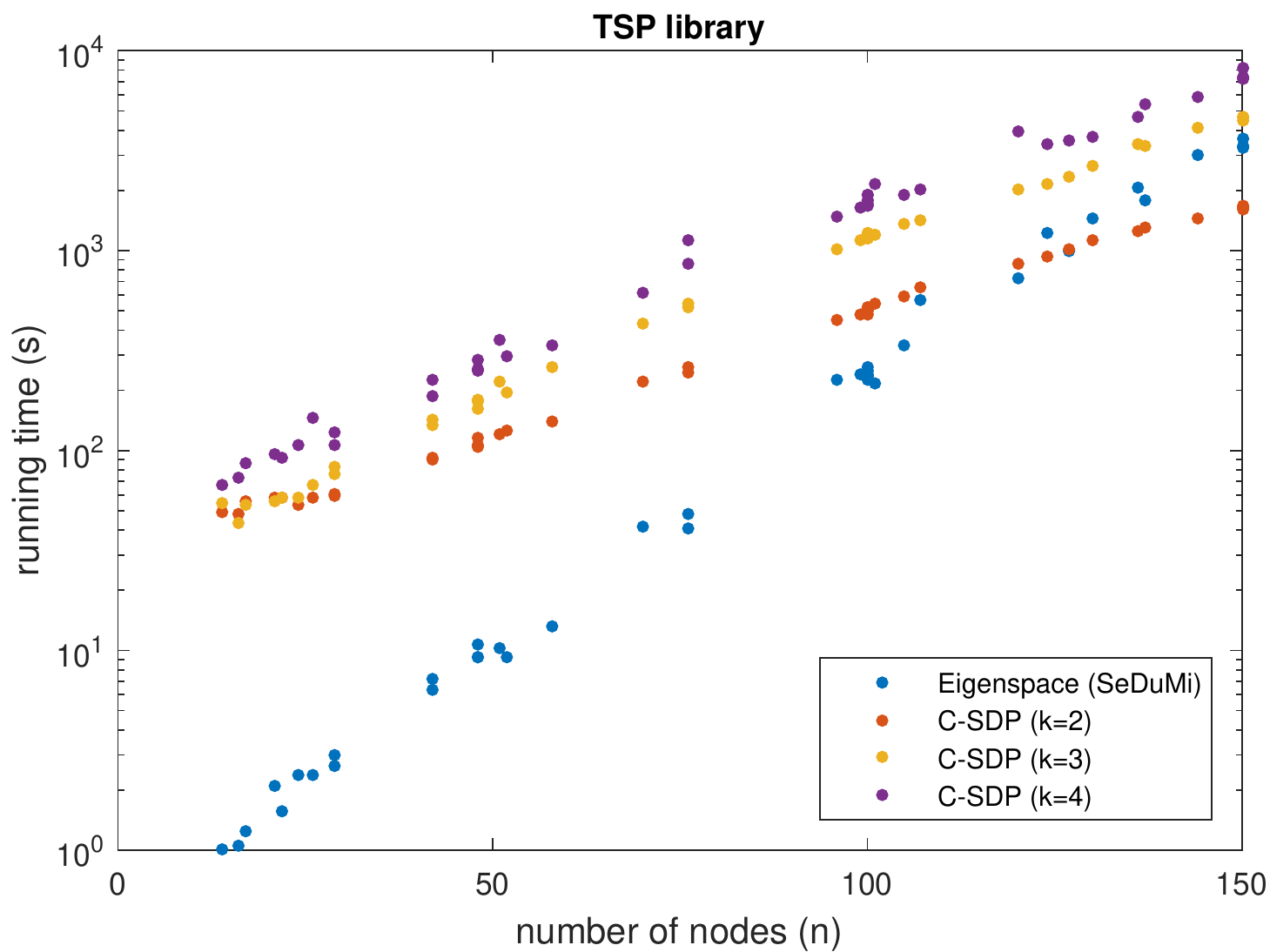}
	\caption{Comparison of run times between C-SDP and Eigenspace (linear program) relexations on problems from the TSP library (with $n\leq 150$). C-SDP instances were ran for 1000 ADMM iterations on 20 processors. Eigenspace instances were solved using SeDuMi \protect\cite{Sedumi}. As the Eigenspace relaxation simplifies to a linear program in the case of TSP, solving the problem using an interior point solver is stil competitive with the ADMM approach used in C-SDP.}
\end{figure}
\begin{figure}
	\centering
	\includegraphics[width=0.75\columnwidth]{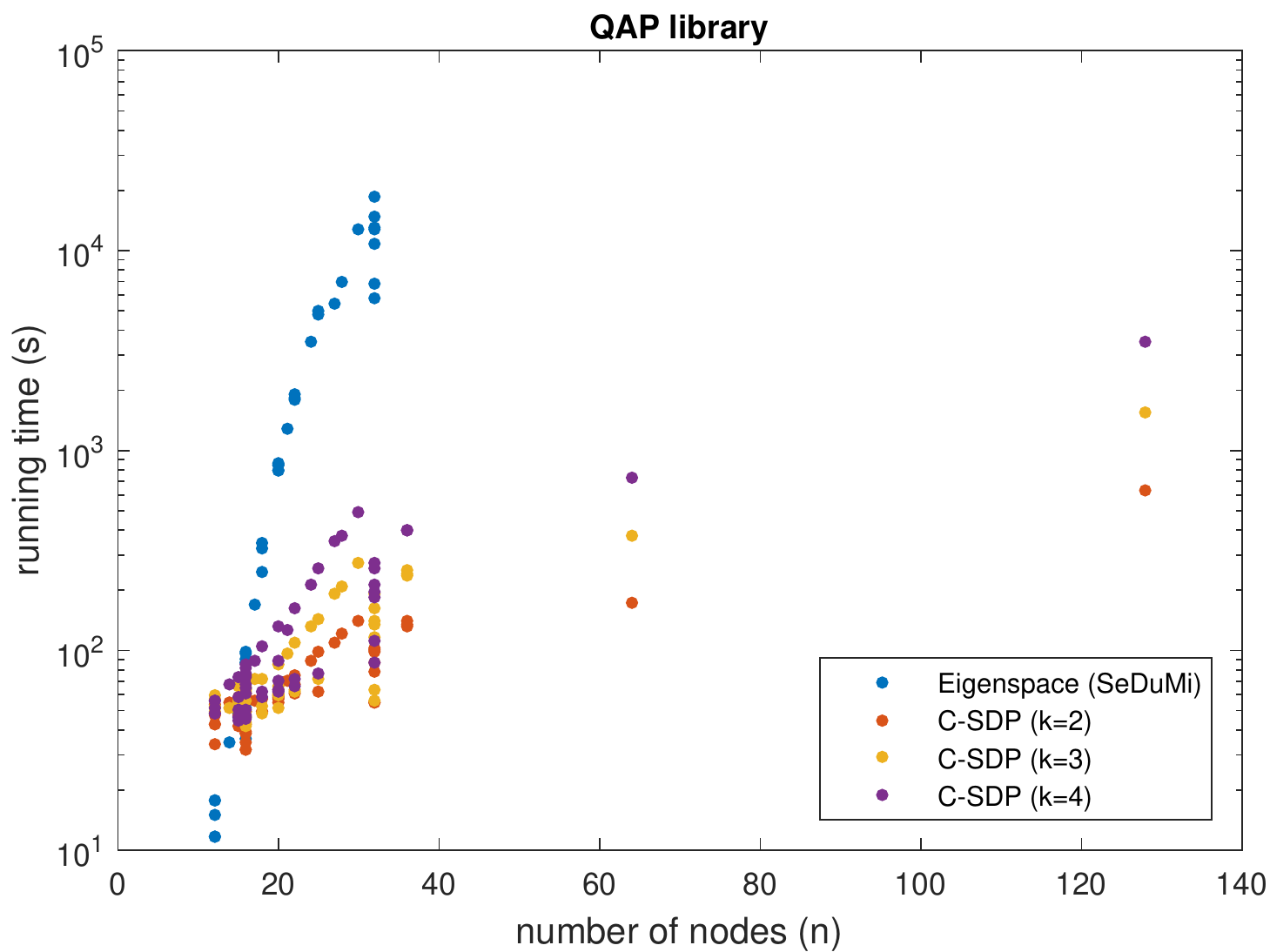}
	\caption{Comparison of run times between C-SDP and Eigenspace relexations on problems from the QAP library (with $n\leq 150$). C-SDP instances were ran for 1000 ADMM iterations on 20 processors. Eigenspace instances were solved using SeDuMi \protect\cite{Sedumi}. In the case of problems from the QAP library, the Eigenspace relaxation no longer simplifies, resulting in longer runtimes. }
\end{figure}

\newpage

\end{document}